\let\chooseClass3   

\ifx\chooseClass1
\IfFileExists{extarticle.cls}{
\documentclass[14pt]{extarticle}
\emergencystretch 7 pt
\setlength{\oddsidemargin}{-17 mm} 
\setlength{\textwidth}{192 mm}     
\setlength{\textheight}{246 mm}    
\setlength{\topmargin}{-31 mm}     
}{
\documentclass[12pt]{article}
\textheight = 8.5in
\textwidth 6.3in
\setlength{\oddsidemargin}{0mm}
\setlength{\topmargin}{-20 mm}
}{}

\makeatletter
\def\@seccntformat#1{\csname the#1\endcsname.\quad}
\renewcommand\section{\@startsection {section}{1}{\z@}%
                                   {-3.5ex \@plus -1ex \@minus -.2ex}%
                                   {2.3ex \@plus.2ex}%
                                   {\normalfont\large\bfseries}}
\renewcommand\subsection{\@startsection{subsection}{2}{\z@}%
                        {3.25ex plus 1ex minus .2ex}{-.5em}%
                        {\normalfont\normalsize\bfseries}}
\renewcommand\subsubsection{\@startsection{subsubsection}{3}{\z@}%
                        {3.25ex plus 1ex minus .2ex}{-.5em}%
                        {\normalfont\normalsize\itseries}}
\@addtoreset{equation}{section}
\makeatother

\usepackage{amsthm}
\newtheoremstyle{boldhead}
{\topsep}
{\topsep}
{\slshape}
{}
{\bfseries}
{.}
{ }
{\thmname{#1}\thmnumber{ #2}\thmnote{ (#3)}}

\newtheoremstyle{boldremark}
{\topsep}
{\topsep}
{\upshape}
{}
{\bfseries}
{.}
{ }
{\thmname{#1}\thmnumber{ #2}\thmnote{ (#3)}}

\swapnumbers

\theoremstyle{boldhead}
\newtheorem{theorem}[subsection]{Theorem}
\newtheorem{corollary}[subsection]{Corollary}
\newtheorem{lemma}[subsection]{Lemma}
\newtheorem{proposition}[subsection]{Proposition}

\theoremstyle{boldremark}
\newtheorem{definition}[subsection]{Definition}
\newtheorem{example}[subsection]{Example}
\newtheorem{examples}[subsection]{Examples}
\newtheorem{remark}[subsection]{Remark}
\newtheorem{problem}[subsection]{Problem}
\newtheorem*{acknowledgement}{Acknowledgements}
\newtheorem*{notation}{Notation}

\fi

\ifx\chooseClass2
\documentclass{tac}
\fi

\ifx\chooseClass3
\documentclass[12pt]{article}
\makeatletter
\def\@seccntformat#1{\csname the#1\endcsname.\quad}
\renewcommand\section{\@startsection {section}{1}{\z@}%
                                   {-3.5ex \@plus -1ex \@minus -.2ex}%
                                   {2.3ex \@plus.2ex}%
                                   {\normalfont\large\bfseries}}
\renewcommand\subsection{\@startsection{subsection}{2}{\z@}%
                        {3.25ex plus 1ex minus .2ex}{-.5em}%
                        {\normalfont\normalsize\bfseries}}
\renewcommand\subsubsection{\@startsection{subsubsection}{3}{\z@}%
                        {3.25ex plus 1ex minus .2ex}{-.5em}%
                        {\normalfont\normalsize\itseries}}
\@addtoreset{equation}{section}
\makeatother

\usepackage{amsthm}
\newtheoremstyle{boldhead}
{\topsep}
{\topsep}
{\slshape}
{}
{\bfseries}
{.}
{ }
{\thmname{#1}\thmnumber{ #2}\thmnote{ (#3)}}

\newtheoremstyle{boldremark}
{\topsep}
{\topsep}
{\upshape}
{}
{\bfseries}
{.}
{ }
{\thmname{#1}\thmnumber{ #2}\thmnote{ (#3)}}

\swapnumbers

\theoremstyle{boldhead}
\newtheorem{theorem}[subsection]{Theorem}
\newtheorem{corollary}[subsection]{Corollary}
\newtheorem{lemma}[subsection]{Lemma}
\newtheorem{proposition}[subsection]{Proposition}

\theoremstyle{boldremark}
\newtheorem{definition}[subsection]{Definition}

\fi

\usepackage{    amsmath,
                amsfonts,
                amssymb,
}

\numberwithin{equation}{subsection}

\IfFileExists{euscript.sty}{\usepackage[mathcal]{euscript}}{}
\IfFileExists{mathrsfs.sty}{\usepackage{mathrsfs}}{\let\mathscr\mathfrak}

\usepackage{ifpdf}
\ifpdf
    \IfFileExists{hyperref.sty}{\usepackage[pdftex]{hyperref}}{}
\else
    \IfFileExists{hyperref.sty}{\usepackage[hypertex]{hyperref}}{}
\fi

\IfFileExists{url.sty}{\usepackage{url}}{}
\providecommand{\url}[1]{{\tt #1}}

\InputIfFileExists{diagrams.tex}{}{%
\IfFileExists{diagrams.sty}{\usepackage{diagrams}}
{\message{Without diagrams you can not process this file!}}}
\diagramstyle[height=2em,balance,righteqno,PostScript=dvips,nohug]

\usepackage{t-angles} 

\def\rhaha{\raise.24ex\hbox{$\rightharpoonup$}\kern-1em\lower.24ex\hbox{$\rightharpoondown$}}%
\def\lhaha{\raise.24ex\hbox{$\leftharpoonup$}\kern-1em\lower.24ex\hbox{$\leftharpoondown$}}%
\def\dhaha{\downharpoonleft\kern-.22em\downharpoonright\kern.02em}%
\def\uhaha{\upharpoonleft\kern-.22em\upharpoonright\kern.02em}
\newarrowhead{twoharpoons}\rhaha\lhaha\dhaha\uhaha

\newarrow{Id}===={twoharpoons}
\newarrow{Epi}----{triangle}
\newarrow{MapsTo}{mapsto}---{->}
\newarrow{Mono}{boldhook}---{->}
\newarrow{TTo}----{->}
\newarrow{Twoar}===={=>}

\newcommand\ZZ{{\mathbb Z}}

\newcommand{\ca}{{\mathcal A}}
\newcommand{\cb}{{\mathcal B}}
\newcommand{\cc}{{\mathcal C}}
\newcommand{\cd}{{\mathcal D}}
\newcommand{\ce}{{\mathcal E}}
\newcommand{\cf}{{\mathcal F}}

\newcommand{\cq}{{\mathcal Q}}
\newcommand{\ct}{{\mathcal T}}
\newcommand{\cu}{{\mathcal U}}

\newcommand{\fu}{{\mathscr U}}

\newcommand{\bull}{{\scriptscriptstyle\bullet}}
\newcommand{\Com}{\mathsf C}
\newcommand{\tdt}{\tens\dots\tens}
\newcommand{\tree}{{\mathfrak t}}
\newcommand{\uni}{{\mathbf i}}

\newcommand{\sS}[2]{\vphantom{#2}#1 #2}
\newcommand{\n}[1]{\nobreakdash-\hspace{0pt}}
\newcommand{\ainf}[1]{$A_\infty$\nobreakdash-\hspace{0pt}}
\newcommand{\ainfu}[1]{$A_\infty^u$\nobreakdash-\hspace{0pt}}
\newcommand{\strAi}{\mathrm{strict}A_\infty}
\newcommand{\Cat}{{\mathcal C}at}

\let\kk\Bbbk

\let\eps\varepsilon
\let\epsilon\varepsilon
\let\ge\geqslant
\let\le\leqslant

\let\tens\otimes

\let\wh\widehat

\DeclareMathOperator\Card{Card}

\DeclareMathOperator\Cone{Cone}
\newcommand{\dgQ}{A_1}

\DeclareMathOperator\id{id}

\DeclareMathOperator\im{Im}
\DeclareMathOperator\inj{in}

\DeclareMathOperator\modul{-mod}
\DeclareMathOperator\Ob{Ob}

\DeclareMathOperator\pr{pr}
\DeclareMathOperator\restr{restr}

\newcommand{\corref}[1]{Corollary~\ref{#1}}
\newcommand{\defref}[1]{Definition~\ref{#1}}

\newcommand{\propref}[1]{Proposition~\ref{#1}}

\newcommand{\secref}[1]{Section~\ref{#1}}
\newcommand{\thmref}[1]{Theorem~\ref{#1}}

\begin{document}
\title{Free $A_\infty$-categories}
    \ifx\chooseClass2
\author{Volodymyr Lyubashenko and Oleksandr Manzyuk}
\address{Institute of Mathematics,
National Academy of Sciences of Ukraine, \\
3 Tereshchenkivska st.,
Kyiv-4, 01601 MSP, Ukraine \\[5pt]
Fachbereich Mathematik,
Postfach 3049,
67653 Kaiserslautern,
Germany
}
\copyrightyear{2006}
\eaddress{lub@imath.kiev.ua\CR manzyuk@mathematik.uni-kl.de}
    \else
\author{Volodymyr Lyubashenko%
\thanks{Institute of Mathematics,
National Academy of Sciences of Ukraine,
3 Tereshchenkivska st.,
Kyiv-4, 01601 MSP,
Ukraine;
lub@imath.kiev.ua}
\ and Oleksandr Manzyuk%
\thanks{Fachbereich Mathematik,
Postfach 3049,
67653 Kaiserslautern,
Germany;
manzyuk@mathematik.uni-kl.de}
}
    \fi
\maketitle

\begin{abstract}
For a differential graded $\kk$\n-quiver $\cq$ we define the free
\ainf-category $\cf\cq$ generated by $\cq$. The main result is that
the restriction \ainf-functor $A_\infty(\cf\cq,\ca)\to A_1(\cq,\ca)$ is
an equivalence, where objects of the last \ainf-category are morphisms
of differential graded $\kk$\n-quivers $\cq\to\ca$.
\end{abstract}

\allowdisplaybreaks[1]

\vspace{1.5em}

\ainf-categories defined by Fukaya~\cite{Fukaya:A-infty} and
Kontsevich~\cite{Kontsevich:alg-geom/9411018} are generalizations of
differential graded categories for which the binary composition is
associative only up to a homotopy. They also generalize \ainf-algebras
introduced by Stasheff~\cite[II]{Stasheff:HomAssoc}. \ainf-functors are
the corresponding generalizations of usual functors, see e.g.
\cite{Fukaya:A-infty,math.RA/9910179}. Homomorphisms of \ainf-algebras
(e.g. \cite{Kadeishvili82}) are particular cases of \ainf-functors.
\ainf-transformations are certain coderivations. Examples of such
structures are encountered in studies of mirror symmetry (e.g.
\cite{Kontsevich:alg-geom/9411018,Fukaya:FloerMirror-II}) and in
homological algebra.

For an \ainf-category there is a notion of units up to a homotopy
(homotopy identity morphisms) \cite{Lyu-AinfCat}. Given two
\ainf-categories $\ca$ and $\cb$, one can construct a third
\ainf-category $A_\infty(\ca,\cb)$, whose objects are \ainf-functors
$f:\ca\to\cb$, and morphisms are \ainf-transformations between such
functors (Fukaya~\cite{Fukaya:FloerMirror-II}, Kontsevich and Soibelman
\cite{KonSoi-AinfCat-NCgeom,KonSoi-book},
Le\-f\`evre-Ha\-se\-ga\-wa~\cite{Lefevre-Ainfty-these}, as well as
\cite{Lyu-AinfCat}). This allows to define a 2\n-category, whose
objects are unital \ainf-categories, 1\n-morphisms are unital
\ainf-functors and 2\n-morphisms are equivalence classes of natural
\ainf-transformations \cite{Lyu-AinfCat}. We continue to study this
2\n-category.

The notations and conventions are explained in the first section. We
also describe $A_N$\n-categories, $A_N$\n-functors and
$A_N$\n-transformations -- truncated at $N<\infty$ versions of
\ainf-categories. For instance, $A_1$\n-categories and $A_1$\n-functors
are differential graded $\kk$\n-quivers and their morphisms. However,
$A_1$\n-transformations bring new 2\n-categorical features to the
theory. In particular, for any differential graded $\kk$\n-quiver $\cq$
and any \ainf-category $\ca$ there is an \ainf-category $A_1(\cq,\ca)$,
whose objects are morphisms of differential graded $\kk$\n-quivers
$\cq\to\ca$, and morphisms are $A_1$\n-transformations. We recall the
terminology related to trees in \secref{sec-trees}.

In the second section we define the free \ainf-category $\cf\cq$
generated by a differential graded $\kk$\n-quiver $\cq$. We classify
functors from a free \ainf-category $\cf\cq$ to an arbitrary
\ainf-category $\ca$ in \propref{pro-extension-f1}. In particular, the
restriction map gives a bijection between the set of strict
\ainf-functors $\cf\cq\to\ca$ and the set of morphisms of differential
graded $\kk$\n-quivers $\cq\to(\ca,m_1)$ (\corref{cor-restr-strAinf}).
We classify chain maps into complexes of transformations whose source
is a free \ainf-category in \propref{pro-chain-to-Ainf-FQ-A}.
Description of homotopies between such chain maps is given in
\corref{cor-chain-to-Ainf-FQ-A-null-homotopic}. Assuming in addition
that $\ca$ is unital, we obtain our main result: the restriction
\ainf-functor $\restr:A_\infty(\cf\cq,\ca)\to A_1(\cq,\ca)$ is an
equivalence (\thmref{thm-restr-equivalence}).

In the third section we interpret $A_\infty(\cf\cq,\_)$ and
$A_1(\cq,\_)$ as strict \ainfu-2-functors $A_\infty^u\to A_\infty^u$.
Moreover, we interpret $\restr:A_\infty(\cf\cq,\_)\to A_1(\cq,\_)$ as
an \ainfu-2-equivalence. In this sense the \ainf-category $\cf\cq$
represents the \ainfu-2-functor $A_1(\cq,\_)$. This is the
2\n-categorical meaning of freeness of $\cf\cq$.

\section{Conventions and preliminaries}\label{sec-convent-nota}
We keep the notations and conventions of
\cite{Lyu-AinfCat,LyuOvs-iResAiFn}, sometimes without explicit
mentioning. Some of the conventions are recalled here.

We assume as in \cite{Lyu-AinfCat,LyuOvs-iResAiFn} that most quivers,
\ainf-categories, etc. are small with respect to some universe $\fu$.

The ground ring $\kk\in\fu$ is a unital associative commutative ring. A
$\kk$\n-module means a $\fu$\n-small $\kk$\n-module.

We use the right operators: the composition of two maps (or morphisms)
$f:X\to Y$ and $g:Y\to Z$ is denoted $fg:X\to Z$; a map is written on
elements as $f:x\mapsto xf=(x)f$. However, these conventions are not
used systematically, and $f(x)$ might be used instead.

$\ZZ$\n-graded $\kk$\n-modules are functions
$X:\ZZ\ni d\mapsto X^d\in\kk\modul$. A simple computation shows that
the product $X=\prod_{\iota\in I}X_\iota$ in the category of
$\ZZ$\n-graded $\kk$\n-modules of a family $(X_\iota)_{\iota\in I}$ of
objects $X_\iota:d\mapsto X_\iota^d$ is given by
$X:\ZZ\ni d\mapsto X^d=\prod_{\iota\in I}X_\iota^d$. Everywhere in this
article the product of graded $\kk$\n-modules means the above product.

If $P$ is a $\ZZ$\n-graded $\kk$\n-module, then $sP=P[1]$ denotes the
same $\kk$\n-module with the grading $(sP)^d=P^{d+1}$. The ``identity''
map $P\to sP$ of degree $-1$ is also denoted $s$. The map $s$ commutes
with the components of the differential in an \ainf-category
(\ainf-algebra) in the following sense: $s^{\tens n}b_n=m_ns$.

Let $\Com=\Com(\kk\modul)$ denote the differential graded category of
complexes of $\kk$-modules. Actually, it is a symmetric closed monoidal
category.

The cone of a chain of a chain map \(\alpha:P\to Q\) of complexes of
$\kk$\n-modules is the graded $\kk$\n-module
\(\Cone(\alpha)=Q\oplus P[1]\) with the differential
\((q,ps)d=(qd^Q+p\alpha,psd^{P[1]})=(qd^Q+p\alpha,-pd^Ps)\).

\subsection{$A_N$-categories.}\label{sec-AN-categories}
For a positive integer $N$ we define some $A_N$\n-notions similarly to
the case $N=\infty$. We may say that all data, equations and
constructions for $A_N$\n-case are the same as in \ainf-case (e.g.
\cite{Lyu-AinfCat}), however, taken only up to level $N$.

A \emph{differential graded $\kk$\n-quiver} $\cq$ is the following
data: a $\fu$\n-small set of objects $\Ob\cq$; a chain complex of
$\kk$\n-modules $\cq(X,Y)$ for each pair of objects $X$, $Y$. A
morphism of differential graded $\kk$\n-quivers $f:\cq\to\ca$ is given
by a map $f:\Ob\cq\to\Ob\ca$, $X\mapsto Xf$ and by a chain map
$\cq(X,Y)\to\ca(Xf,Yf)$ for each pair of objects $X$, $Y$ of $\cq$. The
category of differential graded $\kk$\n-quivers is denoted $\dgQ$.

The category of $\fu$\n-small graded $\kk$\n-linear quivers, whose set
of objects is $S$, admits a monoidal structure with the
tensor product $\ca\times\cb\mapsto\ca\tens\cb$,
$(\ca\tens\cb)(X,Y)=\oplus_{Z\in S}\ca(X,Z)\tens_\kk\cb(Z,Y)$. In
particular, we have tensor powers $T^n\ca=\ca^{\tens n}$ of a given
graded $\kk$\n-quiver $\ca$, such that $\Ob T^n\ca=\Ob\ca$. Explicitly,
\[ T^n\ca(X,Y) = \bigoplus_{X=X_0,X_1,\dots,X_n=Y\in\Ob\ca}
\ca(X_0,X_1)\tens_\kk\ca(X_1,X_2)\tens_\kk\dots\tens_\kk\ca(X_{n-1},X_n).
\]
In particular, $T^0\ca(X,Y)=\kk$ if $X=Y$ and vanishes otherwise. The
graded $\kk$\n-quiver $T^{\le N}\ca=\oplus_{n\ge0}^NT^n\ca$ is called
the restricted tensor coalgebra of $\ca$. It is equipped with the cut
comultiplication
\begin{align*}
\Delta:T^{\le N}\ca(X,Y) &\to \bigoplus_{Z\in\Ob\ca}
T^{\le N}\ca(X,Z)\bigotimes_\kk T^{\le N}\ca(Z,Y), \\
h_1\tens h_2\tens\dots\tens h_n &\mapsto \sum_{k=0}^n
h_1\tens\dots\tens h_k\bigotimes h_{k+1}\tens\dots\tens h_n,
\end{align*}
and the counit
 $\eps=\bigl(T^{\le N}\ca(X,Y) \rTTo^{\pr_0} T^0\ca(X,Y)\to\kk\bigr)$,
where the last map is $\id_\kk$ if $X=Y$, or $0$ if $X\ne Y$ (and
$T^0\ca(X,Y)=0$). We write $T\ca$ instead of $T^{\le\infty}\ca$. If
$g:T\ca\to T\cb$ is a map of $\kk$\n-quivers, then $g_{ac}$ denotes its
matrix coefficient
$T^a\ca \rMono^{\inj_a} T\ca \rTTo^g T\cb \rEpi^{\pr_c} T^c\cb$. The
matrix coefficient $g_{a1}$ is abbreviated to $g_a$.

\begin{definition}
An \emph{$A_N$\n-category} $\ca$ consists of the following data: a
graded $\kk$\n-quiver $\ca$; a system of $\kk$\n-linear maps of degree 1
\[b_n:s\ca(X_0,X_1)\tens s\ca(X_1,X_2)\tens\dots\tens s\ca(X_{n-1},X_n)
\to s\ca(X_0,X_n), \qquad 1\le n\le N,
\]
such that for all $1\le k\le N$
\begin{equation}
\sum_{r+n+t=k} (1^{\tens r}\tens b_n\tens1^{\tens t})b_{r+1+t} = 0 :
T^ks\ca \to s\ca.
\label{eq-b-b-0}
\end{equation}
\end{definition}

The system $b_n$ is interpreted as a (1,1)-coderivation
$b:T^{\le N}s\ca\to T^{\le N}s\ca$ of degree 1 determined by
\begin{equation*}
b_{kl}=(b\big|_{T^ks\ca})\pr_l : T^ks\ca \to T^ls\ca, \qquad
b_{kl} = \sum_{\substack{r+n+t=k\\r+1+t=l}}
1^{\tens r}\tens b_n\tens1^{\tens t}, \qquad k,l\le N,
\end{equation*}
which is a differential.

\begin{definition}
A \emph{pointed cocategory homomorphism} consists of the following
data: $A_N$\n-categories $\ca$ and $\cb$, a map $f:\Ob\ca\to\Ob\cb$ and
a system of $\kk$\n-linear maps of degree 0
\[f_n:s\ca(X_0,X_1)\tens s\ca(X_1,X_2)\tens\dots\tens s\ca(X_{n-1},X_n)
\to s\cb(X_0f,X_nf), \qquad 1\le n\le N.
\]
\end{definition}

The above data are equivalent to a cocategory homomorphism
$f:T^{\le N}s\ca\to T^{\le N}s\cb$ of degree 0 such that
\begin{equation}
f_{01} = (f\big|_{T^0s\ca})\pr_1 =0 : T^0s\ca \to T^1s\cb,
\label{eq-f-T0sA-pr1-0}
\end{equation}
(this condition was implicitly assumed in
\cite[Definition~2.4]{Lyu-AinfCat}). The components of $f$ are
\begin{equation}
f_{kl}=(f\big|_{T^ks\ca})\pr_l : T^ks\ca \to T^ls\cb, \qquad
f_{kl} = \sum_{i_1+\dots+i_l=k}
f_{i_1} \tens f_{i_2} \tens\dots\tens f_{i_l},
\label{eq-f-kl-components}
\end{equation}
where $k,l\le N$. Indeed, the claim follows from the following diagram,
commutative for all $l\ge0$:
\begin{diagram}
T^{\le N}s\ca & \rTTo^f & T^{\le N}s\cb & \rTTo^{\pr_l} & T^ls\cb \\
\dTTo<{\Delta^{(l)}} & = & \dTTo<{\Delta^{(l)}} & = & \dEq \\
(T^{\le N}s\ca)^{\tens l} & \rTTo^{f^{\tens l}} &
(T^{\le N}s\cb)^{\tens l} & \rTTo^{\pr_1^{\tens l}} & (s\cb)^{\tens l}
\end{diagram}
where $\Delta^{(0)}=\eps$, $\Delta^{(1)}=\id$, $\Delta^{(2)}=\Delta$
and $\Delta^{(l)}$ means the cut comultiplication, iterated $l-1$
times. Notice that condition \eqref{eq-f-T0sA-pr1-0} can be written as
$f_0=0$.

\begin{definition}
An \emph{$A_N$\n-functor} $f:\ca\to\cb$ is a pointed cocategory
homomorphism, which commutes with the differential $b$, that is, for
all $1\le k\le N$
\begin{equation*}
\sum_{l>0;i_1+\dots+i_l=k}
(f_{i_1} \tens f_{i_2} \tens\dots\tens f_{i_l}) b_l =
\sum_{r+n+t=k} (1^{\tens r}\tens b_n\tens1^{\tens t}) f_{r+1+t} :
T^ks\ca \to s\cb.
\end{equation*}
\end{definition}

We are interested mostly in the case $N=1$. Clearly, $A_1$\n-categories
are differential graded quivers and $A_1$\n-functors are their
morphisms. In the case of one object these reduce to chain complexes
and chain maps. The following notion seems interesting even in this
case.

\begin{definition}
An $A_N$\n-transformation $r:f\to g:\ca\to\cb$ of degree $d$ consists
of the following data: $A_N$\n-categories $\ca$ and $\cb$; pointed
cocategory homomorphisms $f,g:T^{\le N}s\ca\to T^{\le N}s\cb$ (or
$A_N$\n-functors $f,g:\ca\to\cb$); a system of $\kk$\n-linear maps of
degree $d$
\[ r_n:s\ca(X_0,X_1)\tens s\ca(X_1,X_2)\tens\dots\tens s\ca(X_{n-1},X_n)
\to s\cb(X_0f,X_ng), \qquad 0\le n\le N.
\]
\end{definition}

To give a system $r_n$ is equivalent to specifying an
$(f,g)$-coderivation $r:T^{\le N}s\ca\to T^{\le N+1}s\cb$ of degree $d$
\begin{align}
r_{kl} &= (r\big|_{T^ks\ca})\pr_l : T^ks\ca \to T^ls\cb,
\qquad k\le N,\;l\le N+1 \notag \\
r_{kl} &= \sum_{\substack{q+1+t=l\\i_1+\dots+i_q+n+j_1+\dots+j_t=k}}
f_{i_1} \tens\dots\tens f_{i_q} \tens r_n\tens
g_{j_1} \tens\dots\tens g_{j_t},
\label{eq-r-kl-components}
\end{align}
that is, a $\kk$\n-quiver morphism $r$, satisfying
$r\Delta=\Delta(f\tens r+r\tens g)$. This follows from the commutative
diagram
\begin{diagram}
T^{\le N}s\ca & \rTTo^r & T^{\le N+1}s\cb & \rTTo^{\pr_l} & T^ls\cb \\
\dTTo<{\Delta^{(l)}} & = & \dTTo<{\Delta^{(l)}} & = & \dEq \\
(T^{\le N}s\ca)^{\tens l}
& \rTTo^{\sum_{q+1+t=l}f^{\tens q}\tens r\tens g^{\tens t}} &
(T^{\le N+1}s\cb)^{\tens l} & \rTTo^{\pr_1^{\tens l}} &(s\cb)^{\tens l}
\end{diagram}

Let $\ca$, $\cb$ be $A_N$\n-categories, and let
$f^0,f^1,\dots,f^n:T^{\le N}s\ca\to T^{\le N}s\cb$ be pointed
cocategory homomorphisms. Consider coderivations $r_1$, \dots, $r_n$ as
in
\[ f^0 \rTTo^{r^1} f^1 \rTTo^{r^2} \dots f^{n-1} \rTTo^{r^n} f^n :
T^{\le N}s\ca \to T^{\le N}s\cb.
\]
We construct the following system of $\kk$\n-linear maps
$\theta_{kl}:T^ks\ca\to T^ls\cb$, $k\le N$, $l\le N+n$ of degree
$\deg r^1+\dots+\deg r^n$ from these data:
\begin{equation}
\theta_{kl} = \sum f^0_{i^0_1}\tens\dots\tens f^0_{i^0_{m_0}}\tens
r^1_{j_1}\tens f^1_{i^1_1}\tens\dots\tens f^1_{i^1_{m_1}}\tens\dots\tens
r^n_{j_n}\tens f^n_{i^n_1}\tens\dots\tens f^n_{i^n_{m_n}},
\label{eq-theta-kl-sum}
\end{equation}
where summation is taken over all terms with
\[ m_0+m_1+\dots+m_n+n=l, \quad i^0_1+\dots+i^0_{m_0}+j_1+
i^1_1+\dots+i^1_{m_1}+\dots+j_n+i^n_1+\dots+i^n_{m_n} = k.
\]
Equivalently, we write
\begin{equation*}
\theta_{kl} =
\sum_{\substack{m_0+m_1+\dots+m_n+n=l\\p_0+j_1+p_1+\dots+j_n+p_n=k}}
f^0_{p_0m_0}\tens r^1_{j_1}\tens f^1_{p_1m_1}\tens\dots\tens
r^n_{j_n}\tens f^n_{p_nm_n}.
\end{equation*}
The component $\theta_{kl}$ vanishes unless $n\le l\le k+n$. If $n=0$,
then $\theta_{kl}$ is expansion~\eqref{eq-f-kl-components} of $f^0$. If
$n=1$, then $\theta_{kl}$ is expansion~\eqref{eq-r-kl-components} of
$r^1$.

Given an $A_K$\n-category $\ca$ and an $A_{K+N}$\n-category $\cb$,
$1\le K,N\le\infty$, we construct an $A_N$\n-category $A_K(\ca,\cb)$
out of these. The objects of $A_K(\ca,\cb)$ are $A_K$\n-functors
$f:\ca\to\cb$. Given two such functors $f,g:\ca\to\cb$ we define the
graded $\kk$\n-module $A_K(\ca,\cb)(f,g)$ as the space of all
$A_K$\n-transformations $r:f\to g$, namely,
\begin{multline*}
[A_K(\ca,\cb)(f,g)]^{d+1} \\
= \{ r:f\to g \mid \text{ $A_K$\n-transformation }
r:T^{\le K}s\ca\to T^{\le K+1}s\cb \text{ has degree } d \}.
\end{multline*}
The system of differentials $B_n$, $n\le N$, is defined as follows:
\begin{gather*}
B_1: A_K(\ca,\cb)(f,g)\to A_K(\ca,\cb)(f,g),
\quad r\mapsto (r)B_1 = [r,b] = rb-(-)^rbr, \\
[(r)B_1]_k = \sum_{i_1+\dots+i_q+n+j_1+\dots+j_t=k}
(f_{i_1} \tens\dots\tens f_{i_q} \tens r_n\tens
g_{j_1} \tens\dots\tens g_{j_t})b_{q+1+t} \\
-(-)^r \sum_{\alpha+n+\beta=k}
(1^{\tens\alpha}\tens b_n\tens1^{\tens\beta})r_{\alpha+1+\beta},
\qquad k\le K, \\
B_n: A_K(\ca,\cb)(f^0,f^1)\tens\dots\tens A_K(\ca,\cb)(f^{n-1},f^n)\to
A_K(\ca,\cb)(f^0,f^n), \\
r^1\tens\dots\tens r^n\mapsto (r^1\tens\dots\tens r^n)B_n,
\text{ for } 1<n\le N,
\end{gather*}
where the last $A_K$\n-transformation is defined by its components:
\begin{equation*}
[(r^1\tens\dots\tens r^n)B_n]_k =
\sum_{l=n}^{n+k} (r^1\tens\dots\tens r^n)\theta_{kl}b_l, \qquad k\le K.
\end{equation*}

The category of graded $\kk$\n-linear quivers admits a symmetric
monoidal structure with the tensor product
$\ca\times\cb\mapsto\ca\boxtimes\cb$, where
$\Ob\ca\boxtimes\cb=\Ob\ca\times\Ob\cb$ and
$(\ca\boxtimes\cb)\bigl((X,U),(Y,V)\bigr)=\ca(X,Y)\tens_\kk\cb(U,V)$.
The same tensor product was denoted $\tens$ in \cite{Lyu-AinfCat}, but
we will keep notation $\ca\tens\cb$ only for tensor product from
\secref{sec-AN-categories}, defined when $\Ob\ca=\Ob\cb$. The two
tensor products obey

\textbf{Distributivity law.} Let $\ca$, $\cb$, $\cc$, $\cd$ be graded
$\kk$\n-linear quivers, such that $\Ob\ca=\Ob\cb$ and $\Ob\cc=\Ob\cd$.
Then the middle four interchange map $1\tens c\tens1$ is an isomorphism
of quivers
\begin{equation}
(\ca\tens\cb)\boxtimes(\cc\tens\cd) \rTTo^\sim
(\ca\boxtimes\cc)\tens(\cb\boxtimes\cd),
\label{eq-Distributivity-law}
\end{equation}
identity on objects.

Indeed, the both quivers in \eqref{eq-Distributivity-law} have the same
set of objects $R\times S$, where $R=\Ob\ca=\Ob\cb$ and
$S=\Ob\cc=\Ob\cd$. Let $X,Z\in R$ and $U,W\in S$. The sets of morphisms
from $(X,U)$ to $(Z,W)$ are isomorphic via
\begin{diagram}
\bigl((\ca\tens\cb)\boxtimes(\cc\tens\cd)\bigr)\bigl((X,U),(Z,W)\bigr)
= \\
\bigl(\oplus_{Y\in R}\ca(X,Y)\tens_\kk\cb(Y,Z)\bigr)\tens_\kk
\bigl(\oplus_{V\in S}\cc(U,V)\tens_\kk\cd(V,W)\bigr) \\
\dTTo>\wr \\
\oplus_{(Y,V)\in R\times S}\ca(X,Y)\tens_\kk\cb(Y,Z)\tens_\kk
\cc(U,V)\tens_\kk\cd(V,W) \\
\dTTo>{1\tens c\tens1} \\
\oplus_{(Y,V)\in R\times S}\ca(X,Y)\tens_\kk\cc(U,V)\tens_\kk
\cb(Y,Z)\tens_\kk\cd(V,W) \\
= \bigl((\ca\boxtimes\cc)\tens(\cb\boxtimes\cd)\bigr)
\bigl((X,U),(Z,W)\bigr).
\end{diagram}

The notion of a pointed cocategory homomorphism extends to the case of
several arguments, that is, to degree 0 cocategory homomorphisms
 $\psi:T^{\le L^1}s\cc^1\boxtimes\dots\boxtimes T^{\le L^q}s\cc^q\to
 T^{\le N}s\cb$,
where $N\ge L^1+\dots+L^q$. We always assume that
$\psi_{00\dots0}:T^0s\cc^1\boxtimes\dots\boxtimes T^0s\cc^q\to s\cb$ vanishes.
We call $\psi$ an $A$\n-functor if it commutes with the differential,
that is,
\[ (b\boxtimes1\boxtimes\dots\boxtimes1
+ 1\boxtimes b\boxtimes\dots\boxtimes1 +\dots
+ 1\boxtimes1\boxtimes\dots\boxtimes b)\psi=\psi b.
\]
For example, the map
 $\alpha:T^{\le K}s\ca\boxtimes T^{\le N}sA_K(\ca,\cb)\to
 T^{\le K+N}s\cb$,
 $a\boxtimes r^1\tens\dots\tens r^n\mapsto
 a.[(r^1\tens\dots\tens r^n)\theta]$,
is an $A$\n-functor.

\begin{proposition}[cf. Proposition~5.5 of \cite{Lyu-AinfCat}]
 \label{prop-phi-C1-Cq-psi-ainf}
Let $\ca$ be an $A_K$\n-category, let $\cc^t$ be an
$A_{L^t}$\n-category for $1\le t\le q$, and let $\cb$ be an
$A_N$\n-category, where $N\ge K+L^1+\dots+L^q$. For any $A$\n-functor
 $\phi:T^{\le K}s\ca\boxtimes T^{\le L^1}s\cc^1\boxtimes\dots\boxtimes
 T^{\le L^q}s\cc^q\to T^{\le N}s\cb$
there is a unique $A$\n-functor
 $\psi:T^{\le L^1}s\cc^1\boxtimes\dots\boxtimes T^{\le L^q}s\cc^q
 \to T^{\le N-K}sA_K(\ca,\cb)$,
such that
\begin{equation*}
\phi = \bigl(
 T^{\le K}s\ca\boxtimes T^{\le L^1}s\cc^1\boxtimes\dots\boxtimes T^{\le L^q}s\cc^q
\rTTo^{1\boxtimes\psi} T^{\le K}s\ca\boxtimes T^{\le N-K}sA_K(\ca,\cb)
\rTTo^\alpha T^{\le N}s\cb\bigr).
\end{equation*}
\end{proposition}

Let $\ca$ be an $A_N$\n-category, let $\cb$ be an $A_{N+K}$\n-category,
and let $\cc$ be an $A_{N+K+L}$\n-category. The above proposition
implies the existence of an $A$\n-functor (cf.
\cite[Proposition~4.1]{Lyu-AinfCat})
\[ M: T^{\le K}sA_N(\ca,\cb)\boxtimes T^{\le L}sA_{N+K}(\cb,\cc)
\to T^{\le K+L}sA_N(\ca,\cc),
\]
in particular, \((1\boxtimes B+B\boxtimes1)M=MB\). It has the components
\begin{equation*}
M_{nm}=M\big|_{T^n\boxtimes T^m}\pr_1:
T^nsA_N(\ca,\cb)\boxtimes T^msA_{N+K}(\cb,\cc) \to sA_N(\ca,\cc),
\end{equation*}
$n\le K$, $m\le L$. We have $M_{00}=0$ and $M_{nm}=0$ for $m>1$. If
$m=0$ and $n$ is positive, $M_{n0}$ is given by the formula:
\begin{align*}
M_{n0}:sA_N(\ca,\cb)(f^0,f^1)\tens\dots\tens sA_N(\ca,\cb)(f^{n-1},f^n)
\boxtimes\kk_{g^0} &\to sA_N(\ca,\cc)(f^0g^0,f^ng^0), \\
r^1\tens\dots\tens r^n\boxtimes1 &\mapsto
(r^1\tens\dots\tens r^n\mid g^0)M_{n0},
\end{align*}
\[ [(r^1\tens\dots\tens r^n\mid g^0)M_{n0}]_k =
\sum_{l=n}^{n+k}(r^1\tens\dots\tens r^n)\theta_{kl}g^0_l,\qquad k\le N,
\]
where $\mid$ separates the arguments in place of $\boxtimes$. If $m=1$,
then $M_{n1}$ is given by the formula:
\begin{multline*}
M_{n1}:sA_N(\ca,\cb)(f^0,f^1)\tens\dots\tens sA_N(\ca,\cb)(f^{n-1},f^n)
\boxtimes sA_{N+K}(\cb,\cc)(g^0,g^1) \\
\to sA_N(\ca,\cc)(f^0g^0,f^ng^1),\qquad r^1\tens\dots\tens r^n\boxtimes t^1
\mapsto (r^1\tens\dots\tens r^n\boxtimes t^1)M_{n1},
\end{multline*}
\[ [(r^1\tens\dots\tens r^n\boxtimes t^1)M_{n1}]_k =
\sum_{l=n}^{n+k}(r^1\tens\dots\tens r^n)\theta_{kl}t^1_l,\qquad k\le N.
\]
Note that equations
\begin{equation*}
[(r^1\tens\dots\tens r^n)B_n]_k =
[(r^1\tens\dots\tens r^n\boxtimes b)M_{n1}]_k
-(-)^{r^1+\dots+r^n} [(b\boxtimes r^1\tens\dots\tens r^n)M_{1n}]_k
\end{equation*}
imply that
\begin{align*}
(r^1\tens\dots\tens r^n)B_n &= (r^1\tens\dots\tens r^n\boxtimes b)M_{n1}
-(-)^{r^1+\dots+r^n} (b\boxtimes r^1\tens\dots\tens r^n)M_{1n}, \\
B &= (1\boxtimes b)M - (b\boxtimes1)M: \id\to\id: A_N(\ca,\cb)\to A_N(\ca,\cb).
\end{align*}

\propref{prop-phi-C1-Cq-psi-ainf} implies the existence of a unique
$A_L$\n-functor
\[ A_N(\ca,\_): A_{N+K}(\cb,\cc) \to A_K(A_N(\ca,\cb),A_N(\ca,\cc)),
\]
such that
\begin{multline*}
M = \bigl[ T^{\le K}sA_N(\ca,\cb)\boxtimes T^{\le L}sA_{N+K}(\cb,\cc)
\rTTo^{1\boxtimes A_N(\ca,\_)} \\
T^{\le K}sA_N(\ca,\cb)\boxtimes T^{\le L}sA_K(A_N(\ca,\cb),A_N(\ca,\cc))
\rTTo^\alpha T^{\le K+L}sA_N(\ca,\cc)\bigr].
\end{multline*}
The $A_L$\n-functor $A_N(\ca,\_)$ is strict, cf.
\cite[Proposition~6.2]{Lyu-AinfCat}.

Let $\ca$ be an $A_N$\n-category, and let $\cb$ be a unital
\ainf-category with a unit transformation $\uni^\cb$. Then
$A_N(\ca,\cb)$ is a unital \ainf-category with the unit transformation
$(1\boxtimes\uni^\cb)M$ (cf. \cite[Proposition~7.7]{Lyu-AinfCat}). The unit
element for an object $f\in\Ob A_N(\ca,\cb)$ is
$\sS{_f}\uni^{A_N(\ca,\cb)}_0:\kk\to(sA_N)^{-1}(\ca,\cb)$,
$1\mapsto f\uni^\cb$.

When $\ca$ is an $A_K$\n-category and $N<K$, we may forget part of its
structure and view $\ca$ as an $A_N$\n-category. If furthermore, $\cb$
is an $A_{K+L}$\n-category, we have the restriction strict
$A_L$\n-functor $\restr_{K,N}:A_K(\ca,\cb)\to A_N(\ca,\cb)$. To prove
the results mentioned above, we notice that they are restrictions of
their \ainf-analogs to finite $N$. Since the proofs of \ainf-results
are obtained in \cite{Lyu-AinfCat} by induction, an inspection shows
that the proofs of the above $A_N$\n-statements are obtained as a
byproduct.

\subsection{Trees.}\label{sec-trees}
Since the notions related to trees might be interpreted with some
variations, we give precise definitions and fix notation. A \emph{tree}
is a non-empty connected graph without cycles. A vertex which belongs
to only one edge is called \emph{external}, other vertices are
\emph{internal}. A \emph{plane tree} is a tree equipped for each
internal vertex $v$ with a cyclic ordering of the set $E_v$ of edges,
adjacent to $v$. Plane trees can be drawn on an oriented plane in a
unique way (up to an ambient isotopy) so that the cyclic ordering of
each $E_v$ agrees with the orientation of the plane. An external vertex
distinct from the root is called \emph{input vertex}.

A \emph{rooted tree} is a tree with a distinguished external vertex,
called \emph{root}. The set of vertices $V(t)$ of a rooted tree $t$ has
a canonical ordering: $x\preccurlyeq y$ iff the minimal path connecting
the root with $y$ contains $x$. A \emph{linearly ordered tree} is a
rooted tree $t$ equipped with a linear order $\le$ of the set of
internal vertices $IV(t)$, such that for all internal vertices $x$, $y$
the relation $x\preccurlyeq y$ implies $x\le y$. For each vertex
$v\in V(t)-\{\text{root}\}$ of a rooted tree, the set $E_v$ has a
distinguished element $e_v$ -- the beginning of a minimal path from $v$
to the root. Therefore, for each vertex $v\in V(t)-\{\text{root}\}$ of
a rooted plane tree, the set $E_v$ admits a unique linear order $<$,
for which $e_v$ is minimal and the induced cyclic order is the given
one. An internal vertex $v$ has degree $d$, if $\Card(E_v)=d+1$.

For any $y\in V(t)$ let $P_y=\{x\in V(t)\mid x\preccurlyeq y\}$. With
each plane rooted tree $t$ is associated a linearly ordered tree
$t_<=(t,\le)$ as follows. If $x,y\in IV(t)$ are such that
$x\not\preccurlyeq y$ and $y\not\preccurlyeq x$, then $P_x\cap P_y=P_z$
for a unique $z\in IV(t)$, distinct from $x$ and $y$. Let
$a\in E_z-\{e_z\}$ (resp. $b\in E_z-\{e_z\}$) be the beginning of the
minimal path connecting $z$ and $x$ (resp. $y$). If $a<b$, we set
$x<y$. Graphically we $<$-order the internal vertices by height. Thus,
an internal vertex $x$ on the left is depicted lower than a
$\preccurlyeq$\n-incomparable internal vertex $y$ on the right:
\[
\hstretch 150
\begin{tangles}{rl}
\vstretch 33
\begin{tangles}[b]{l}
\nodel{x}\node \\
\se1\node \\
\step\se1\node\nodeld{a}\noder{z}
\end{tangles}
&
\begin{tangles}[b]{r}
\node\noder{y} \\
\hh\node\noded{b}\sw1 \\
\sw1\step
\end{tangles}
\\
\n\noderu{e_z} \\
\hh\n \\
\hh\id \\
\object{\text{root}}
\end{tangles}
\quad,\qquad a<b \implies x<y.
\]

A \emph{forest} is a sequence of plane rooted trees. Concatenation of
forests is denoted $\sqcup$. The vertical composition $F_1\cdot F_2$ of
forests $F_1$, $F_2$ is well-defined if the sum of lengths of sequences
$F_1$ and $F_2$ equals the number of external vertices of $F_2$. These
operations allow to construct any tree from elementary ones
\[ 1=\;
\begin{tangle}
\id
\end{tangle}
\quad, \qquad \text{and} \qquad \tree_k=
\vstretch 50
\begin{tangle}
\nw2\nw1\node\Put(1,17)[cb]{\scriptscriptstyle\cdots}\ne1\ne2\\
\step[2]\id
\end{tangle}
\qquad (k\text{ input vertices}).
\]
Namely, any linearly ordered tree $(t,\le)$ has a unique presentation
of the form
\begin{equation}
(t,\le)=(1^{\sqcup\alpha_1}\sqcup\tree_{k_1}\sqcup1^{\sqcup\beta_1})\cdot
(1^{\sqcup\alpha_2}\sqcup\tree_{k_2}\sqcup1^{\sqcup\beta_2})\cdot
\ldots\cdot \tree_{k_N},
\label{eq-ord-tree-decomp-forest}
\end{equation}
where $N=|t|\overset{\text{def}}=\Card(IV(t))$ is the number of
internal vertices. Here
\[ 1^{\sqcup\alpha}\sqcup\tree_k\sqcup1^{\sqcup\beta} =\;
\begin{tangles}[b]{ccc}
\object{\overbrace{\step[4]}^\alpha} & \object{\overbrace{\step[4]}^k}
& \object{\overbrace{\step[4]}^\beta} \\
\id\step\id\step\nodeu{\dots}\step\id\step\id & \step
\vstretch 50
\begin{tangles}[b]{c}
\nw2\nw1\node\Put(1,17)[cb]{\scriptscriptstyle\cdots}\ne1\ne2\\
\id
\end{tangles}
\step & \id\step\id\step\nodeu{\dots}\step\id\step\id
\end{tangles}
\quad.
\]
In \eqref{eq-ord-tree-decomp-forest} the highest vertex is indexed by
1, the lowest -- by $N$.

\section{\texorpdfstring{Properties of free $A_\infty$-categories}
 {Properties of free A-infinity categories}}
\subsection{\texorpdfstring{Construction of a free $A_\infty$-category.}
 {Construction of a free A-infinity category.}}
The category $\strAi$ has \ainf-categories as objects and strict
\ainf-functors as morphisms. There is a functor $\cu:\strAi\to\dgQ$,
$\ca\mapsto(\ca,m_1)$ which sends an \ainf-category to the underlying
differential graded $\kk$\n-quiver, forgetting all higher
multiplications. Following Kontsevich and
Soibelman~\cite{KonSoi-AinfCat-NCgeom} we are going to prove that $\cu$
has a left adjoint functor $\cf:\dgQ\to\strAi$, $\cq\mapsto\cf\cq$.
The \ainf-category $\cf\cq$ is called free. Below we describe its
structure for an arbitrary differential graded $\kk$\n-quiver $\cq$. We
shall work with its shift $(s\cq,d)$.

Let us define an \ainf-category $\cf\cq$ via the following data. The
class of objects $\Ob\cf\cq$ is $\Ob\cq$. The $\ZZ$\n-graded
$\kk$\n-modules of morphisms between $X,Y\in\Ob\cq$ are
\begin{align*}
s\cf\cq(X,Y) &= \oplus_{n\ge1} \oplus_{t\in\ct^n_{\ge2}}
s\cf_t\cq(X,Y), \\
s\cf_t\cq(X,Y) &= \oplus_{X_0,\dots,X_n\in\Ob\cq}^{X_0=X,\,X_n=Y}
s\cq(X_0,X_1)\tens\dots\tens s\cq(X_{n-1},X_n)\bigl[-|t|\bigr],
\end{align*}
where $\ct^n_{\ge2}$ is the class of plane rooted trees with $n+1$
external vertices, such that $\Card(E_v)\ge3$ for all $v\in IV(t)$. We
use the following convention: if $M$, $N$ are (differential) graded
$\kk$\n-modules, then%
\footnote{Another gauge choice $(M\tens N)[1]=M[1]\tens N$, $s=s\tens1$
seems less convenient.}
\begin{gather*}
(M\tens N)[k] = M\tens\bigl(N[k]\bigr), \\
\bigl(M\tens N \rTTo^{s^k} (M\tens N)[k]\bigr) =
\bigl(M\tens N \rTTo^{1\tens s^k} M\tens(N[k])\bigr).
\end{gather*}
The quiver $\cf\cq$ is equipped with the following operations. For
$k>1$ the operation $b_k$ is a direct sum of maps
\begin{equation}
b_k = s^{|t_1|}\tens\dots\tens s^{|t_{k-1}|}\tens s^{|t_k|-|t|}:
s\cf_{t_1}\cq(Y_0,Y_1)\tens\dots\tens s\cf_{t_k}\cq(Y_{k-1},Y_k)
\to s\cf_t\cq(Y_0,Y_k),
\label{eq-bk-free-Acat-ssssss}
\end{equation}
where $t=(t_1\sqcup\dots\sqcup t_k)\cdot\tree_k$. In particular,
$|t|=|t_1|+\dots+|t_k|+1$. The operation $b_1$ restricted to
$s\cf_t\cq$ is
\begin{equation}
b_1 = d\oplus(-1)^{\beta(t')}s^{-1}: s\cf_t\cq(X,Y) \to
s\cf_t\cq(X,Y) \oplus \bigoplus_{t'=t+\text{edge}} s\cf_{t'}\cq(X,Y),
\label{eq-b1-free-Acat}
\end{equation}
where the sum extends over all trees $t'\in\ct^n_{\ge2}$ with a
distinguished edge $e$, such that contracting $e$ we get $t$ from $t'$.
The sign is determined by
\[ \beta(t') = \beta(t',e) = 1 + h(\text{highest vertex of }e), \]
where an isomorphism of ordered sets
\[ h: IV(t'_<) \rTTo^\sim \bigl[1,|t'|\bigr]\cap\ZZ
\]
is simply the height of a vertex in the linearly ordered tree $t'_<$,
canonically associated with $t'$. In \eqref{eq-b1-free-Acat} $d$ means
$d\tens1\tens\dots\tens1+\dots+1\tens\dots\tens d\tens1
+1\tens\dots\tens1\tens d$,
where the last $d$ is
$d_{s\cq[-|t|]}=(-)^{|t|}s^{|t|}\cdot d_{s\cq}\cdot s^{-|t|}$, as
usual. According to our conventions, $s^{-1}$ in
\eqref{eq-b1-free-Acat} means $1^{\tens n-1}\tens s^{-1}$.

\begin{proposition}
$\cf\cq$ is an \ainf-category.
\end{proposition}

\begin{proof}
First we prove that $b_1^2=0$ on $s\cf_t\cq$. Indeed,
\begin{multline*}
b_1^2 = d^2\oplus (-)^{\beta(t',e)}(s^{-1}d+ds^{-1}) \oplus
\bigl[(-1)^{\beta(t'_1,e_1)+\beta(t'',e_2)}
+(-1)^{\beta(t'_2,e_2)+\beta(t'',e_1)}\bigr]s^{-2}: \\
s\cf_t\cq \to s\cf_t\cq \oplus \bigoplus_{t'=t+e} s\cf_{t'}\cq
\oplus \bigoplus_{t''=t+e_1+e_2} s\cf_{t''}\cq,
\end{multline*}
where $t''$ contains two distinguished edges $e_1$, $e_2$, contraction
along which gives $t$; $t'_2$ is $t''$ contracted along $e_1$, and
$t'_1$ is $t''$ contracted along $e_2$. We may assume that highest
vertex of $e_1$ is lower than highest vertex of $e_2$ in $t''_<$. Then
$\beta(t'_1,e_1)=\beta(t'',e_1)$ and
$\beta(t'',e_2)=\beta(t'_2,e_2)+1$, hence,
\[ (-1)^{\beta(t'_1,e_1)+\beta(t'',e_2)}
+(-1)^{\beta(t'_2,e_2)+\beta(t'',e_1)}=0.
\]
Obviously, $d^2=0$ and $s^{-1}d+ds^{-1}=0$, hence, $b_1^2=0$.

Let us prove for each $n>1$ that
\begin{multline*}
b_nb_1 + \sum_{p=1}^n(1^{\tens p-1}\tens b_1\tens1^{\tens n-p})b_n
+ \sum_{\substack{\alpha+k+\beta=n\\\alpha+\beta>0}}^{k>1}
(1^{\tens\alpha}\tens b_k\tens1^{\tens\beta})b_{\alpha+1+\beta} =0: \\
s\cf_{t_1}\cq\tens\dots\tens s\cf_{t_n}\cq \to s\cf_t\cq
\oplus \bigoplus_{p,t'_p} s\cf_{t'}\cq \oplus
\bigoplus_{\substack{\alpha+k+\beta=n\\\alpha+\beta>0}}^{k>1,\,t''}
s\cf_{t''}\cq,
\end{multline*}
where
\begin{gather}
t = (t_1\sqcup\dots\sqcup t_n)\cdot\tree_n =\;
\vstretch 50
\hstretch 200
\begin{tanglec}
\ffbox{1}{t_1}\step\ffbox{1}{t_2}\step\ffbox{1}{t_{n-1}}\step\ffbox{1}{t_n} \\
\nw3\step\nw1\node\Put(1,17)[cb]{\scriptscriptstyle\cdots}\ne1\step\ne3 \\
\id
\end{tanglec}
\;, \label{eq-T-T1-Tn-tan} \\
t' = (t_1\sqcup\dots\sqcup t'_p\sqcup\dots\sqcup t_n)\cdot\tree_n =\;
\vstretch 50
\hstretch 200
\begin{tanglec}
\ffbox{1}{t_1}\step\ffbox{1}{t'_p}\step\ffbox{1}{t_n} \\
\nw2\nw1\n\ne1\ne2 \\
\id
\end{tanglec}
\;, \notag
\end{gather}
\begin{multline*}
t'' = (t_1\sqcup\dots\sqcup t_n)\cdot
(1^{\sqcup\alpha}\sqcup\tree_k\sqcup1^{\sqcup\beta})
\cdot\tree_{\alpha+1+\beta}
\\
=\;
\def\FillCircDiam{1}
\hstretch 300
\begin{tangles}{llcrr}
\HH\ffbox{1}{t_1}&\ffbox{1}{t_\alpha}&\ffbox{1}{t_{\alpha+1}}\step
\ffbox{1}{t_{\alpha+k}}&\ffbox{1}{t_{\alpha+k+1}}&\ffbox{1}{t_n} \\
\hstep\nw4&\hstep{\hstr{375}\nw2}\step[-0.25]&\hstr{600}\hd\n{\hh\s}\hdd
&\step[-0.25]{\hstr{375}\ne2}\hstep&\ne4\hstep \\
&&\hh\id&&
\end{tangles}
\;,
\end{multline*}
and contraction of $t'_p$ along distinguished edge $e_p$ gives $t_p$.
According to the three types of summands in the target, the required
equation follows from anticommutativity of the following three
diagrams:
\begin{diagram}
s\cf_{t_1}\cq\tens\dots\tens s\cf_{t_n}\cq &
\rTTo^{s^{|t_1|}\tens\dots\tens s^{|t_{n-1}|}\tens s^{|t_n|-|t|}}_{b_n}
& s\cf_t\cq \\
\dTTo<{1^{\tens n-1}\tens d+\dots+d\tens1^{\tens n-1}} & - & \dTTo>d \\
s\cf_{t_1}\cq\tens\dots\tens s\cf_{t_n}\cq &
\rTTo^{s^{|t_1|}\tens\dots\tens s^{|t_{n-1}|}\tens s^{|t_n|-|t|}}_{b_n}
& s\cf_t\cq
\end{diagram}
that is,
\begin{multline*}
(1^{\tens n-1}\tens d+\dots+d\tens1^{\tens n-1})
(s^{|t_1|}\tens\dots\tens s^{|t_{n-1}|}\tens s^{|t_n|-|t|}) \\
+(s^{|t_1|}\tens\dots\tens s^{|t_{n-1}|}\tens s^{|t_n|-|t|})
(1^{\tens n-1}\tens d+\dots+d\tens1^{\tens n-1}) = 0;
\end{multline*}
\begin{diagram}[nobalance]
s\cf_{t_1}\cq\tens\dots\tens s\cf_{t_p}\cq\tens\dots\tens s\cf_{t_n}\cq &
\rTTo^{s^{|t_1|}\tens\dots\tens s^{|t_{n-1}|}\tens s^{|t_n|-|t|}}
& s\cf_t\cq \\
\dTTo>{(-)^{\beta(t'_p)}1^{\tens p-1}\tens s^{-1}\tens1^{\tens n-p}} & -
& \dTTo<{\hspace*{-8em}(-)^{1+|t_1|+\dots+|t_{p-1}|+\beta(t'_p)}s^{-1}} \\
s\cf_{t_1}\cq\tens\dots\tens s\cf_{t'_p}\cq\tens\dots\tens s\cf_{t_n}\cq &
\rTTo_{s^{|t_1|}\tens\dots\tens s^{|t_{p-1}|}\tens s^{|t_p|+1}\tens
s^{|t_{p+1}|}\tens\dots\tens s^{|t_{n-1}|}\tens s^{|t_n|-|t|-1}}
& s\cf_{t'}\cq \\
\end{diagram}
that is,
\begin{multline*}
(-1)^{\beta(t'_p)}(1^{\tens p-1}\tens s^{-1}\tens1^{\tens n-p})
(s^{|t_1|}\tens\dots\tens s^{|t_{p-1}|}\tens s^{|t_p|+1}\tens
s^{|t_{p+1}|}\tens\dots\tens s^{|t_{n-1}|}\tens s^{|t_n|-|t|-1}) \\
+ (s^{|t_1|}\tens\dots\tens s^{|t_{n-1}|}\tens s^{|t_n|-|t|})
(-1)^{1+|t_1|+\dots+|t_{p-1}|+\beta(t'_p)}(1^{\tens n-1}\tens s^{-1}) =0,
\end{multline*}
in the particular case $p=n$ it holds as well;
\begin{diagram}[nobalance,objectstyle=\scriptstyle]
s\cf_{t_1}\cq\tens\dots\tens s\cf_{t_n}\cq &
\rTTo^{s^{|t_1|}\tens\dots\tens s^{|t_{n-1}|}\tens s^{|t_n|-|t|}}
& s\cf_t\cq \\
\dTTo<{1^{\tens\alpha}\tens s^{|t_{\alpha+1}|}\tens\dots}>{\tens
s^{|t_{\alpha+k-1}|}\tens s^{-|t_{\alpha+1}|-\dots-|t_{\alpha+k-1}|-1}
\tens1^{\tens\beta}\hspace*{-8em}} & - &
\dTTo<{\hspace*{-5em}(-)^{1+|t_1|+\dots+|t_\alpha|}s^{-1}} \\
s\cf_{t_1}\cq\tens\dots\tens s\cf_{t_\alpha}\cq\tens s\cf_{\hat t}\cq
\tens s\cf_{t_{\alpha+k+1}}\cq\tens\dots\tens s\cf_{t_n}\cq
& \rTTo_{s^{|t_1|}\tens\dots\tens s^{|t_\alpha|}\tens
s^{|\hat t|}\tens s^{|t_{\alpha+k+1}|}
\tens\dots\tens s^{|t_{n-1}|}\tens s^{|t_n|-|t|-1}}
& s\cf_{t''}\cq \\
\end{diagram}
where $\hat t=(t_{\alpha+1}\sqcup\dots\sqcup t_{\alpha+k})\cdot\tree_k$,
that is,
\begin{multline*}
(1^{\tens\alpha}\tens s^{|t_{\alpha+1}|}\tens\dots\tens
s^{|t_{\alpha+k-1}|}\tens s^{-|t_{\alpha+1}|-\dots-|t_{\alpha+k-1}|-1}
\tens1^{\tens\beta})\cdot \\
(s^{|t_1|}\tens\dots\tens s^{|t_\alpha|}\tens1^{\tens k-1}\tens
s^{|t_{\alpha+1}|+\dots+|t_{\alpha+k}|+1}\tens s^{|t_{\alpha+k+1}|}
\tens\dots\tens s^{|t_{n-1}|}\tens s^{|t_n|-|t|-1}) \\
+ (s^{|t_1|}\tens\dots\tens s^{|t_{n-1}|}\tens s^{|t_n|-|t|})
(-1)^{1+|t_1|+\dots+|t_\alpha|}(1^{\tens n-1}\tens s^{-1}) =0,
\end{multline*}
in the particular case $\beta=0$ it holds as well.

Therefore, $\cf\cq$ is an \ainf-category.
\end{proof}

Let us establish a property of free \ainf-categories, which explains
why they are called free.

\begin{proposition}[$A_\infty$-functors from a free $A_\infty$-category]
\label{pro-extension-f1}
Let $\cq$ be a differential graded quiver, and let $\ca$ be an
\ainf-category. Let $f_1:s\cq\to(s\ca,b_1)$ be a chain morphism of
differential graded quivers with the underlying mapping of objects
$\Ob f:\Ob\cq\to\Ob\ca$. Suppose given $\kk$\n-quiver morphisms
$f_k:T^ks\cf\cq\to s\ca$ of degree 0 with the same underlying map
$\Ob f$ for all $k>1$. Then there exists a unique extension of $f_1$ to
a quiver morphism $f_1:s\cf\cq\to s\ca$ such that $(f_1,f_2,\dots)$
are components of an \ainf-functor $f:\cf\cq\to\ca$.
\end{proposition}

\begin{proof}
For each $n>1$ we have to satisfy the equation
\begin{equation}
b_nf_1 = \sum_{i_1+\dots+i_l=n}(f_{i_1}\tdt f_{i_l})b_l
- \sum_{\alpha+k+\beta=n}^{\alpha+\beta>0}
(1^{\tens\alpha}\tens b_k\tens1^{\tens\beta})f_{\alpha+1+\beta}:
T^ns\cf\cq\to s\ca.
\label{eq-bnf1-def-f1}
\end{equation}
It is used to define recursively $f_1$ on $s\cf\cq$. Suppose that
$t_1$, \dots, $t_n$ are trees, $n>1$, and $f_1:s\cf_{t_i}\cq\to s\ca$
is already defined for all $1\le i\le n$. Since
\begin{equation*}
b_n = s^{|t_1|}\tdt s^{|t_{n-1}|}\tens s^{|t_n|-|t|}:
s\cf_{t_1}\cq\tdt s\cf_{t_n}\cq \to s\cf_t\cq
\end{equation*}
is invertible for $t=(t_1\sqcup\dots\sqcup t_n)\cdot\tree_n$,
formula~\eqref{eq-bnf1-def-f1} determines $f_1:s\cf_t\cq\to s\ca$
uniquely as
\[ f_1 = \bigl( s\cf_t\cq \rTTo^{b_n^{-1}}
s\cf_{t_1}\cq\tdt s\cf_{t_n}\cq
\rTTo^{\sum(f_{i_1}\tdt f_{i_l})b_l
- \sum_{\alpha+k+\beta=n}^{\alpha+\beta>0}
(1^{\tens\alpha}\tens b_k\tens1^{\tens\beta})f_{\alpha+1+\beta}}
s\ca \bigr).
\]
This proves uniqueness of the extension of $f_1$.

Let us prove that the cocategory homomorphism $f$ with so defined
components $(f_1,f_2,\dots)$ is an \ainf-functor.
Equations~\eqref{eq-bnf1-def-f1} are satisfied by construction of
$f_1$. So it remains to prove that $f_1$ is a chain map. Equation
$f_1b_1=b_1f_1$ holds on $s\cf_|\cq$ by assumption. We are going to
prove by induction on $|t|$ that it holds on $s\cf_t\cq$. Considering
$t=(t_1\sqcup\dots\sqcup t_n)\cdot\tree_n$, $n>1$, we assume that
$f_1b_1=b_1f_1:s\cf_{t'}\cq\to s\ca$ for all trees $t'$ with
$|t'|<|t|$. To prove that $f_1b_1=b_1f_1:s\cf_t\cq\to s\ca$ it suffices
to show that $b_nf_1b_1=b_nb_1f_1$ for all $n>1$ due to invertibility
of $b_n$. Using \eqref{eq-bnf1-def-f1} and the equation $b^2\pr_1=0$ we
find
\begin{align*}
b_nf_1b_1 - b_nb_1f_1 &=
\sum_{i_1+\dots+i_l=n}(f_{i_1}\tdt f_{i_l})b_lb_1
- \sum_{\alpha+k+\beta=n}^{\alpha+\beta>0}
(1^{\tens\alpha}\tens b_k\tens1^{\tens\beta})f_{\alpha+1+\beta}b_1 \\
&\qquad + \sum_{\alpha+k+\beta=n}^{\alpha+\beta>0}
(1^{\tens\alpha}\tens b_k\tens1^{\tens\beta})b_{\alpha+1+\beta}f_1 \\
&= - \sum_{i_1+\dots+i_l=n}(f_{i_1}\tdt f_{i_l})
\sum_{\gamma+p+\delta=l}^{\gamma+\delta>0}
(1^{\tens\gamma}\tens b_p\tens1^{\tens\delta})b_{\gamma+1+\delta} \\
&\qquad + \sum_{\alpha+k+\beta=n}^{\alpha+\beta>0}
(1^{\tens\alpha}\tens b_k\tens1^{\tens\beta}) \Biggl[
\sum_{j_1+\dots+j_r=\alpha+1+\beta}^{r>1}(f_{j_1}\tdt f_{j_r})b_r \\
&\qquad\qquad\qquad\qquad\qquad
- \sum_{\gamma+p+\delta=\alpha+1+\beta}^{\gamma+\delta>0}
(1^{\tens\gamma}\tens b_p\tens1^{\tens\delta})f_{\gamma+1+\delta}
\Biggr] \\
&= \sum_{r>1} \Biggl[ \sum_{\alpha+k+\beta=n}^{\alpha+\beta>0}
(1^{\tens\alpha}\tens b_k\tens1^{\tens\beta})
\sum_{j_1+\dots+j_r=\alpha+1+\beta}f_{j_1}\tdt f_{j_r} \\
&\qquad\qquad - \sum_{i_1+\dots+i_l=n}(f_{i_1}\tdt f_{i_l})
\sum_{\substack{\gamma+p+\delta=l\\ \gamma+1+\delta=r}}
1^{\tens\gamma}\tens b_p\tens1^{\tens\delta} \Biggr]b_r \\
&\qquad - \sum_{r>1} \Biggl[
\sum_{\substack{\alpha+k+\beta=n\\ \alpha+\beta>0}}
(1^{\tens\alpha}\tens b_k\tens1^{\tens\beta})
\sum_{\substack{\gamma+p+\delta=\alpha+1+\beta\\ \gamma+1+\delta=r}}
1^{\tens\gamma}\tens b_p\tens1^{\tens\delta} \Biggr]f_r.
\end{align*}

Let us show that the expressions in square brackets vanish. The first
one is the matrix coefficient
$bf-fb:s\cf_{t_1}\cq\tdt s\cf_{t_n}\cq\to T^rs\ca$. Indeed, for $r>1$
the inequality $r\le j_1+\dots+j_r=\alpha+1+\beta$ automatically
implies that $\alpha+\beta>0$, so this condition can be omitted. Using
the induction hypothesis one can transform the left hand side of
equation
\begin{multline*}
\sum_{\alpha+k+\beta=n} (1^{\tens\alpha}\tens b_k\tens1^{\tens\beta})
\sum_{j_1+\dots+j_r=\alpha+1+\beta}f_{j_1}\tdt f_{j_r} \\
= \sum_{i_1+\dots+i_l=n}(f_{i_1}\tdt f_{i_l})
\sum_{\substack{\gamma+p+\delta=l\\ \gamma+1+\delta=r}}
1^{\tens\gamma}\tens b_p\tens1^{\tens\delta}:
s\cf_{t_1}\cq\tdt s\cf_{t_n}\cq \to T^rs\ca
\end{multline*}
into the right hand side for all $n,r\ge1$.

The second expression
\begin{equation}
\sum_{\substack{\alpha+k+\beta=n\\ \alpha+\beta>0}}
(1^{\tens\alpha}\tens b_k\tens1^{\tens\beta})
\sum_{\substack{\gamma+p+\delta=\alpha+1+\beta\\ \gamma+1+\delta=r}}
1^{\tens\gamma}\tens b_p\tens1^{\tens\delta}
\label{eq-1b1-1b1-0}
\end{equation}
is the matrix coefficient
\[ (b-b\pr_1)b\pr_r: T^ns\cf\cq \to T^rs\cf\cq
\]
of the endomorphism $(b-b\pr_1)b:Ts\cf\cq\to Ts\cf\cq$. However,
\[ (b-b\pr_1)b\pr_r = b^2\pr_r - b\pr_1b\pr_r = - b\pr_1b\pr_1\pr_r = 0
\]
for $r>1$, because $\pr_1b=\pr_1b\pr_1$. Therefore,
\eqref{eq-1b1-1b1-0} vanishes and equation $b_nf_1b_1=b_nb_1f_1$ is
proven.
\end{proof}

Let $\strAi(\cf\cq,\ca)\subset A_\infty(\cf\cq,\ca)$ be a full
\ainf-subcategory, whose objects are strict \ainf-functors. Recall that
$\Ob\dgQ(\cq,\ca)$ is the set of chain morphisms $\cq\to\ca$ of
differential graded quivers.

\begin{corollary}\label{cor-restr-strAinf}
A chain morphism $f:\cq\to\ca$ admits a unique extension to a strict
\ainf-functor $\widehat{f}:\cf\cq\to\ca$. The maps
$f\mapsto\widehat{f}$ and
\begin{equation*}
\mathrm{restr}: \Ob\strAi(\cf\cq,\ca) \to \Ob\dgQ(\cq,\ca),
\qquad g \mapsto (\Ob g,g_1\big|_{s\cq})
\end{equation*}
are inverse to each other.
\end{corollary}

Indeed, strict \ainf-functors $g$ are distinguished by conditions
$g_k=0$ for $k>1$.

We may view $\strAi$ as a category, whose objects are \ainf-categories
and morphisms are strict \ainf-functors. We may also view $\dgQ$ as a
category consisting of differential graded quivers and their morphisms.
There is a functor $\cu:\strAi\to\dgQ$, $\ca\mapsto(\ca,m_1)$, which
sends an \ainf-category to the underlying differential graded
$\kk$\n-quiver, forgetting all higher multiplications. The restriction
map
\begin{equation}
\mathrm{restr}: \strAi(\cf\cq,\ca) \to \dgQ(\cq,\cu\ca),
\qquad g \mapsto (\Ob g,g_1\big|_{s\cq})
\label{eq-restr-strAinf-dgQ}
\end{equation}
is functorial in $\ca$.

\begin{corollary}
There is a functor $\cf:\dgQ\to\strAi$, $\cq\mapsto\cf\cq$, left
adjoint to $\cu$.
\end{corollary}

\subsection{Explicit formula for the constructed strict
 $A_\infty$-functor.}
Let us obtain a more explicit formula for $\wh{f}_1\big|_{s\cf_t\cq}$.
We define $\wh{f}_1$ for $t=(t_1\sqcup\dots\sqcup t_n)\cdot\tree_n$
recursively by a commutative diagram
\begin{diagram}
s\cf_{t_1}\cq\tens\dots\tens s\cf_{t_n}\cq &
\rTTo^{s^{|t_1|}\tens\dots\tens s^{|t_{n-1}|}\tens s^{|t_n|-|t|}}_{b_n}
& s\cf_t\cq \\
\dTTo<{\wh{f}_1^{\tens n}} & \overset{\text{def}}=:
& \dTTo>{\wh{f}_1} \\
s\ca\tens\dots\tens s\ca & \rTTo^{b_n} & s\ca
\end{diagram}
Notice that the top map is invertible. Here $n$, $t_1$, \dots, $t_n$
are uniquely determined by decomposition~\eqref{eq-T-T1-Tn-tan} of $t$.

Let $(t,\le)$ be a linearly ordered tree with the underlying given
plane rooted tree $t$. Decompose $(t,\le)$ into a vertical composition
of forests as in \eqref{eq-ord-tree-decomp-forest}. Then the following
diagram commutes
\begin{diagram}
s\cq^{\tens n} & \rTTo^{1^{\tens\alpha_1}\tens b_{k_1}\tens1^{\tens\beta_1}}
& s\cq^{\tens\alpha_1}\tens s\cf_{\tree_{k_1}}\cq\tens s\cq^{\tens\beta_1}
& \rTTo^{1^{\tens\alpha_2}\tens b_{k_2}\tens1^{\tens\beta_2}} & \dots &
\rTTo^{b_{k_N}} & s\cf_t\cq \\
\dTTo<{f_1^{\tens n}} && \dTTo<{\hat{f}_1^{\tens\alpha_1+1+\beta_1}}
&& \dTTo<{\hspace*{-4em} \hat{f}_1^{\tens\alpha_2+1+\beta_2}}
\dots \dTTo>{\hat{f}_1^{\tens k_N}} && \dTTo>{\hat{f}_1} \\
s\ca^{\tens n} & \rTTo^{1^{\tens\alpha_1}\tens b_{k_1}\tens1^{\tens\beta_1}}
& s\ca^{\tens\alpha_1}\tens s\ca\tens s\ca^{\tens\beta_1} &
\rTTo^{1^{\tens\alpha_2}\tens b_{k_2}\tens1^{\tens\beta_2}} & \dots &
\rTTo^{b_{k_N}} & s\ca
\end{diagram}
The upper row consists of invertible maps. One can prove by induction
that the composition of maps in the upper row equals $\pm s^{-|t|}$.
When $(t,\le)=t_<$ is the linearly ordered tree, canonically associated
with $t$, then the composition of maps in the upper row equals
$s^{-|t|}$. This is also proved by induction: if $t$ is presented as
$t=(t_1\sqcup\dots\sqcup t_k)\cdot\tree_k$, then the composition of
maps in the upper row is
\begin{multline*}
(1^{\tens k-1}\tens s^{-|t_k|})
(1^{\tens k-2}\tens s^{-|t_{k-1}|}\tens1)\dots
(s^{-|t_1|}\tens1^{\tens k-1})
(s^{|t_1|}\tens\dots\tens s^{|t_{k-1}|}\tens s^{|t_k|-|t|}) \\
= 1^{\tens k-1}\tens s^{-|t|} = s^{-|t|}.
\end{multline*}
Therefore, for an arbitrary tree $t\in\ct^n_{\ge2}$ the map
$\wh{f}_1\big|_{s\cf_t\cq}$ is
\[ \wh{f}_1 =\bigl(s\cf_t\cq \rTTo^{s^{|t|}} s\cq^{\tens n}
\rTTo^{f_1^{\tens n}} s\ca^{\tens n}
\rTTo^{1^{\tens\alpha_1}\tens b_{k_1}\tens1^{\tens\beta_1}}
s\ca^{\tens\alpha_1+1+\beta_1}
\rTTo^{1^{\tens\alpha_2}\tens b_{k_2}\tens1^{\tens\beta_2}} \dots
\rTTo^{b_{k_N}} s\ca \bigr),
\]
where the factors correspond to
decomposition~\eqref{eq-ord-tree-decomp-forest} of $t_<$.

\subsection{Transformations between functors from a free
$A_\infty$-category.}
 Let $\cq$ be a differential graded quiver, and let $\ca$ be an
\ainf-category. Then $A_1(\cq,\ca)$ is an \ainf-category as well. The
differential graded quiver $(sA_1(\cq,\ca),B_1)$ is described as
follows. Objects are chain quiver maps $\phi:(s\cq,b_1)\to(s\ca,b_1)$,
the graded $\kk$\n-module of morphisms $\phi\to\psi$ is the product of
graded $\kk$\n-modules
\[ sA_1(\cq,\ca)(\phi,\psi) = \prod_{X\in\Ob\cq}s\ca(X\phi,X\psi)\times
\prod_{X,Y\in\Ob\cq}\Com\bigl(s\cq(X,Y),s\ca(X\phi,Y\psi)\bigr),
\quad r = (r_0,r_1).
\]
The differential $B_1$ is given by
\begin{align}
(rB_1)_0 &= r_0b_1, \notag \\
(rB_1)_1 &=
r_1b_1 + (\phi_1\tens r_0)b_2 + (r_0\tens\psi_1)b_2 -(-)^rb_1r_1.
\label{eq-rB10-rB11}
\end{align}
Restrictions $\phi,\psi:\cq\to\ca$ of arbitrary \ainf-functors
$\phi,\psi:\cf\cq\to\ca$ to $\cq$ are $A_1$\n-functors (chain quiver
maps).

\begin{proposition}\label{pro-chain-to-Ainf-FQ-A}
Let $\phi,\psi:\cf\cq\to\ca$ be \ainf-functors. For an arbitrary
complex $P$ of $\kk$\n-modules chain maps
$u:P\to sA_\infty(\cf\cq,\ca)(\phi,\psi)$ are in bijection with the
following data: $(u',u_k)_{k>1}$
\begin{enumerate}
\item a chain map $u':P\to sA_1(\cq,\ca)(\phi,\psi)$,

\item $\kk$\n-linear maps
\[ u_k: P \to \prod_{X,Y\in\Ob\cq}
\Com\bigl((s\cf\cq)^{\tens k}(X,Y),s\ca(X\phi,Y\psi)\bigr)
\]
of degree 0 for all $k>1$.
\end{enumerate}
The bijection maps $u$ to $(u',u_k)_{k>1}$, where $u_k=u\cdot\pr_k$ and
\begin{equation}
u' = \bigl(P \rTTo^u sA_\infty(\cf\cq,\ca)(\phi,\psi) \rEpi^\restr
sA_1(\cf\cq,\ca)(\phi,\psi) \rEpi^\restr sA_1(\cq,\ca)(\phi,\psi)
\bigr).
\label{eq-u'-urestr2}
\end{equation}
The inverse bijection can be recovered from the recurrent formula
\begin{multline*}
(-)^pb^{\cf\cq}_k(pu_1) = -(pd)u_k + \sum_{a+q+c=k}^{\alpha,\beta}
(\phi_{a\alpha}\tens pu_q\tens\psi_{c\beta})b^\ca_{\alpha+1+\beta} \\
-(-)^p \sum_{\alpha+q+\beta=k}^{\alpha+\beta>0} (1^{\tens\alpha}\tens
b^{\cf\cq}_q\tens1^{\tens\beta})(pu_{\alpha+1+\beta}):
(s\cf\cq)^{\tens k} \to s\ca,
\end{multline*}
where $k>1$, \(p\in P\), and \(\phi_{a\alpha}\), \(\psi_{c\beta}\) are
matrix elements of $\phi$, $\psi$.
\end{proposition}

\begin{proof}
Since the $\kk$\n-module of $(\phi,\psi)$-coderivations
$sA_\infty(\cf\cq,\ca)(\phi,\psi)$ is a product, $\kk$\n-linear maps
$u:P\to sA_\infty(\cf\cq,\ca)(\phi,\psi)$ of degree 0 are in bijection
with sequences of $\kk$\n-linear maps $(u_k)_{k\ge0}$ of degree 0:
\begin{alignat*}2
u_0: P &\to \prod_{X\in\Ob\cq}s\ca(X\phi,X\psi), & p &\mapsto pu_0, \\
u_k: P &\to \prod_{X,Y\in\Ob\cq}
\Com\bigl((s\cf\cq)^{\tens k}(X,Y),s\ca(X\phi,Y\psi)\bigr),
&\qquad p &\mapsto pu_k,
\end{alignat*}
for $k\ge1$. The complex
$\Phi_0=(sA_\infty(\cf\cq,\ca)(\phi,\psi),B_1)$ admits a filtration by
subcomplexes
\[ \Phi_n = 0\times\dots\times0\times \prod_{k=n}^\infty
\prod_{X,Y\in\Ob\cq}
\Com\bigl((s\cf\cq)^{\tens k}(X,Y),s\ca(X\phi,Y\psi)\bigr).
\]
In particular, $\Phi_2$ is a subcomplex, and
\[ \Phi_0/\Phi_2 = \prod_{X\in\Ob\cq}s\ca(X\phi,X\psi) \times
\prod_{X,Y\in\Ob\cq}\Com\bigl(s\cf\cq(X,Y),s\ca(X\phi,Y\psi)\bigr)
\]
is the quotient complex with differential~\eqref{eq-rB10-rB11}. Since
$(s\cf\cq,b_1)$ splits into a direct sum of two subcomplexes
$s\cq\oplus\bigl(\oplus_{|t|>0}s\cf_t\cq\bigr)$, the complex
$\Phi_0/\Phi_2$ has a subcomplex
\[ \bigl(0\times \prod_{X,Y\in\Ob\cq}
\Com\bigl(\oplus_{|t|>0}s\cf_t\cq(X,Y),s\ca(X\phi,Y\psi)\bigr),
[\_,b_1]\bigr).
\]
The corresponding quotient complex is $sA_1(\cq,\ca)(\phi,\psi)$. The
resulting quotient map
$\restr_1:sA_\infty(\cf\cq,\ca)(\phi,\psi)\to sA_1(\cq,\ca)(\phi,\psi)$
is the restriction map. Denoting $u'=u\cdot\restr_1$, we get the
discussed assignment $u\mapsto(u',u_n)_{n>1}$. The claim is that if $u$
is a chain map, then the missing part
\[ u''_1: P \to \prod_{X,Y\in\Ob\cq}
\Com\bigl(\oplus_{|t|>0}s\cf_t\cq(X,Y),s\ca(X\phi,Y\psi)\bigr),
\]
of $u_1=u'_1\times u''_1$ is recovered in a unique way.

Let us prove that the map $u\mapsto(u',u_n)_{n>1}$ is injective. The
chain map $u$ satisfies $pdu=puB_1$ for all $p\in P$. That is,
$pdu_k=(puB_1)_k$ for all $k\ge0$. Since
$puB_1=(pu)b^\ca-(-)^pb^{\cf\cq}(pu)$, these conditions can be
rewritten as
\begin{equation}
pdu_k = \sum_{a+q+c=k}^{\alpha,\beta}
(\phi_{a\alpha}\tens pu_q\tens\psi_{c\beta})b^\ca_{\alpha+1+\beta}
-(-)^p \sum_{\alpha+q+\beta=k}
(1^{\tens\alpha}\tens b_q\tens1^{\tens\beta})(pu_{\alpha+1+\beta}),
\label{eq-pduk-pub-bpu}
\end{equation}
where $\phi_{a\alpha}:T^as\cf\cq(X,Y)\to T^\alpha s\ca(X\phi,Y\phi)$
are matrix elements of $\phi$, and $\psi_{c\beta}$ are matrix elements
of $\psi$. The same formula can be rewritten as
\begin{multline}
(-)^pb^{\cf\cq}_k(pu_1) = -(pd)u_k + \sum_{a+q+c=k}^{\alpha,\beta}
(\phi_{a\alpha}\tens pu_q\tens\psi_{c\beta})b^\ca_{\alpha+1+\beta} \\
-(-)^p \sum_{\alpha+q+\beta=k}^{\alpha+\beta>0} (1^{\tens\alpha}\tens
b^{\cf\cq}_q\tens1^{\tens\beta})(pu_{\alpha+1+\beta}):
s\cf_{t_1}\cq\tdt s\cf_{t_k}\cq \to s\ca.
\label{eq-bFQk-pu1}
\end{multline}
When $k>1$, the map
$b^{\cf\cq}_k:s\cf_{t_1}\cq\tdt s\cf_{t_k}\cq\to s\cf_t\cq$,
$t=(t_1\sqcup\dots\sqcup t_k)\tree_k$ is invertible, thus,
$pu_1:s\cf_t\cq\to s\ca$ in the left hand side is determined in a
unique way by $u_0$, $u_n$ for $n>1$ and by
$pu_1:s\cf_{t_i}\cq\to s\ca$, $1\le i\le k$, occurring in the right
hand side. Since the restriction $u_1'$ of $u_1$ to $s\cf_|\cq=s\cq$ is
known by 1), the map $u_1''$ is recursively recovered from
$(u_0,u_1',u_n)_{n>1}$.

Let us prove that the map $u\mapsto(u',u_n)_{n>1}$ is surjective. Given
$(u_0,u_1',u_n)_{n>1}$ we define maps $u_1''$ of degree 0 recursively
by \eqref{eq-bFQk-pu1}. This implies equation~\eqref{eq-pduk-pub-bpu}
for $k>1$. For $k=0$ this equation in the form $pdu_0=pu_0b_1$ holds
due to condition 1). It remains to prove
equation~\eqref{eq-pduk-pub-bpu} for $k=1$:
\begin{multline}
(pd)u_1 = (pu_1)b^\ca_1 + (\phi_1\tens pu_0)b^\ca_2
+ (pu_0\tens\psi_1)b^\ca_2 -(-)^pb_1(pu_1): \\
s\cf_t\cq(X,Y) \to s\ca(X\phi,Y\psi)
\label{eq-pdu1-pu1bA1}
\end{multline}
for all trees $t\in\ct_{\ge2}$. For $t=|$ it holds due to assumption
1). Let $N>1$ be an integer. Assume that
equation~\eqref{eq-pdu1-pu1bA1} holds for all trees $t\in\ct_{\ge2}$
with the number of input leaves $\inj(t)<N$. Let $t\in\ct^N_{\ge2}$ be
a tree (with $\inj(t)=N$). Then $t=(t_1\sqcup\dots\sqcup t_k)\tree_k$
for some $k>1$ and some trees $t_i\in\ct_{\ge2}$, $\inj(t_i)<N$. For
such $t$ equation~\eqref{eq-pdu1-pu1bA1} is equivalent to
\begin{multline*}
(-)^pb_k(pd)u_1 = (-)^pb_k(pu_1)b^\ca_1
+(-)^pb_k(\phi_1\tens pu_0)b^\ca_2 +(-)^pb_k(pu_0\tens\psi_1)b^\ca_2 \\
+ \sum_{\gamma+j+\delta=k}^{\gamma+\delta>0}
(1^{\tens\gamma}\tens b_j\tens1^{\tens\delta})b_{\gamma+1+\delta}(pu_1)
: s\cf_{t_1}\cq\tdt s\cf_{t_k}\cq \to s\ca.
\end{multline*}
Substituting definition~\eqref{eq-bFQk-pu1} of $u_1$ we turn the above
equation into an identity
\begin{align}
& - \sum_{a+q+c=k}^{\alpha,\beta}
(\phi_{a\alpha}\tens pdu_q\tens\psi_{c\beta})b^\ca_{\alpha+1+\beta}
\label{eq-phiaalpha-511} \\
&-(-)^p \sum_{\alpha+q+\beta=k}^{\alpha+\beta>0}
(1^{\tens\alpha}\tens b_q\tens1^{\tens\beta})(pdu_{\alpha+1+\beta})
\label{eq-1bq1pdu-512} \\
&= -(pdu_k)b_1 \label{eq-pdukb1-513} \\
&+ \sum_{a+q+c=k}^{\alpha,\beta}
(\phi_{a\alpha}\tens pu_q\tens\psi_{c\beta})b_{\alpha+1+\beta}b_1
\notag \\
&-(-)^p \sum_{\alpha+q+\beta=k}^{\alpha+\beta>0}
(1^{\tens\alpha}\tens b_q\tens1^{\tens\beta})(pu_{\alpha+1+\beta})b_1
\label{eq-1bq1pub1-51*1} \\
&+(-)^p b_k(\phi_1\tens pu_0)b_2 +(-)^p b_k(pu_0\tens\psi_1)b_2
\label{eq-bfpub-bpuwb-51*2} \\
&+ (-)^p \sum_{\gamma+j+\delta=k}^{\gamma+\delta>0}
(1^{\tens\gamma}\tens b_j\tens1^{\tens\delta})
\biggl[ -pdu_{\gamma+1+\delta} \label{eq-1br1pdu-514} \\
&\qquad\qquad+ \sum_{a+q+c=\gamma+1+\delta}^{\alpha,\beta}
(\phi_{a\alpha}\tens pu_q\tens\psi_{c\beta})b_{\alpha+1+\beta}
\label{eq-fpuwb-51*3} \\
&\qquad\qquad-(-)^p \sum_{\alpha+q+\beta=\gamma+1+\delta}^{\alpha+\beta>0}
(1^{\tens\alpha}\tens b_q\tens1^{\tens\beta})(pu_{\alpha+1+\beta})
\biggr], \label{eq-1bq1pu-515}
\end{align}
whose validity we are going to prove now. First of all, terms
\eqref{eq-1bq1pdu-512} and \eqref{eq-1br1pdu-514} cancel each other.
Term~\eqref{eq-1bq1pu-515} vanishes because for an arbitrary integer
$g$ the sum
\begin{equation}
\sum_{\substack{\gamma+j+\delta=k\\
\alpha+q+\beta=\gamma+1+\delta}}^{\alpha+1+\beta=g}
(1^{\tens\gamma}\tens b_j\tens1^{\tens\delta})
(1^{\tens\alpha}\tens b_q\tens1^{\tens\beta})
\label{eq-1b1-1b1-vanishes}
\end{equation}
is the matrix coefficient $b^2=0:T^ks\cf\cq\to T^gs\cf\cq$, thus, it
vanishes. Notice that condition $\alpha+\beta>0$ in
\eqref{eq-1bq1pu-515} automatically implies $\gamma+\delta>0$.
Furthermore, term~\eqref{eq-pdukb1-513} cancels one of the terms of
sum~\eqref{eq-phiaalpha-511}. In the remaining terms of
\eqref{eq-phiaalpha-511} we may use the induction assumptions and
replace $pdu_q$ with the right hand side of \eqref{eq-pduk-pub-bpu}. We
also absorb terms~\eqref{eq-bfpub-bpuwb-51*2} into
sum~\eqref{eq-fpuwb-51*3}, allowing $\gamma=\delta=0$ in it and
allowing simultaneously $\alpha=\beta=0$ in \eqref{eq-1bq1pub1-51*1} to
compensate for the missing term $b_k(pu_1)b_1$:
\begin{align}
&- \sum_{a+q+c=k}^{\alpha+\beta>0}\sum_{e+j+f=q}^{\gamma,\delta}
\bigl[\phi_{a\alpha}\tens(\phi_{e\gamma}\tens pu_j\tens\psi_{f\delta})
b^\ca_{\gamma+1+\delta}\tens\psi_{c\beta}\bigr]b^\ca_{\alpha+1+\beta}
\label{eq-phifpuwbwb-521} \\
&+(-)^p \sum_{a+q+c=k}^{\alpha+\beta>0} \sum_{\gamma+j+\delta=q}
\bigl[\phi_{a\alpha}\tens(1^{\tens\gamma}\tens b_j\tens1^{\tens\delta})
(pu_{\gamma+1+\delta})\tens\psi_{c\beta}\bigr]b^\ca_{\alpha+1+\beta}
\label{eq-phi1b1puwb-522} \\
&= \sum_{a+q+c=k}^{\alpha,\beta}
(\phi_{a\alpha}\tens pu_q\tens\psi_{c\beta})
b^\ca_{\alpha+1+\beta}b^\ca_1 \label{eq-fpuwbb-523} \\
&-(-)^p \sum_{\alpha+q+\beta=k}
(1^{\tens\alpha}\tens b_q\tens1^{\tens\beta})(pu_{\alpha+1+\beta})
b^\ca_1 \label{eq-1bq1pubA1-524} \\
&+(-)^p \sum_{\gamma+j+\delta=k}
\sum_{a+q+c=\gamma+1+\delta}^{\alpha,\beta}
(1^{\tens\gamma}\tens b_j\tens1^{\tens\delta})
(\phi_{a\alpha}\tens pu_q\tens\psi_{c\beta})b^\ca_{\alpha+1+\beta}.
\notag
\end{align}
Recall that $\phi_{a0}$ vanish for all $a$ except $a=0$. Therefore, we
may absorb term~\eqref{eq-fpuwbb-523} into
sum~\eqref{eq-phifpuwbwb-521} and term~\eqref{eq-1bq1pubA1-524} into
sum~\eqref{eq-phi1b1puwb-522}, allowing terms with $\alpha=\beta=0$ in
them. Denote $r=pu\in sA_\infty(\cf\cq,\ca)(\phi,\psi)$. The
proposition follows immediately form the following

\begin{lemma}\label{lem-ffrpp-identity}
For all $r\in sA_\infty(\cf\cq,\ca)(\phi,\psi)$ and all $k\ge0$ we have
\begin{align}
&- \sum_{a+q+c=k}^{\alpha,\beta}\sum_{e+j+f=q}^{\gamma,\delta}
\bigl[\phi_{a\alpha}\tens(\phi_{e\gamma}\tens r_j\tens\psi_{f\delta})
b^\ca_{\gamma+1+\delta}\tens\psi_{c\beta}\bigr]b^\ca_{\alpha+1+\beta}
\notag \\
&+(-)^r \sum_{a+q+c=k}^{\alpha,\beta} \sum_{\gamma+j+\delta=q}
\bigl[\phi_{a\alpha}\tens(1^{\tens\gamma}\tens b_j\tens1^{\tens\delta})
r_{\gamma+1+\delta}\tens\psi_{c\beta}\bigr]b^\ca_{\alpha+1+\beta}
\notag \\
&= (-)^r \sum_{\gamma+j+\delta=k}
\sum_{a+q+c=\gamma+1+\delta}^{\alpha,\beta}
(1^{\tens\gamma}\tens b_j\tens1^{\tens\delta})
(\phi_{a\alpha}\tens r_q\tens\psi_{c\beta})b^\ca_{\alpha+1+\beta}.
\label{eq-1b1frwb-906}
\end{align}
\end{lemma}

\begin{proof}
Sum~\eqref{eq-1b1frwb-906} is split into three sums accordingly to
output of $b_j$ being an input of $\phi_{a\alpha}$ or $r_q$ or
$\psi_{c\beta}$:
\begin{align}
&- \sum_{a+e+j+f+c=k}^{\alpha,\beta,\gamma,\delta}
(\phi_{a\alpha}\tens\phi_{e\gamma}\tens r_j\tens\psi_{f\delta}
\tens\psi_{c\beta})
(1^{\tens\alpha}\tens b^\ca_{\gamma+1+\delta}\tens1^{\tens\beta})
b^\ca_{\alpha+1+\beta} \label{eq-ffpuww1b1b-531} \\
&+(-)^r \sum_{a+\gamma+j+\delta+c=k}^{\alpha,\beta}
(1^{\tens a+\gamma}\tens b_j\tens1^{\tens\delta+c})
(\phi_{a\alpha}\tens r_{\gamma+1+\delta}\tens\psi_{c\beta})
b^\ca_{\alpha+1+\beta} \label{eq-1b1fpuwb-532} \\
&= (-)^r \sum_{x+q+c=k}^{a,\alpha,\beta}
(b_{xa}\phi_{a\alpha}\tens r_q\tens\psi_{c\beta})
b^\ca_{\alpha+1+\beta} \label{eq-bfpuwb-533} \\
&+(-)^r \sum_{a+y+c=k}^{\alpha,q,\beta}
(\phi_{a\alpha}\tens b_{yq}r_q\tens\psi_{c\beta})
b^\ca_{\alpha+1+\beta} \label{eq-fbpuwb-534} \\
&+ \sum_{a+q+z=k}^{\alpha,\beta,c}
(\phi_{a\alpha}\tens r_q\tens b_{zc}\psi_{c\beta})
b^\ca_{\alpha+1+\beta}. \label{eq-fpubwb-535}
\end{align}
Here $b_{xa}:T^xs\cf\cq\to T^as\cf\cq$ is a matrix element of
$b^{\cf\cq}$. Terms~\eqref{eq-1b1fpuwb-532} and \eqref{eq-fbpuwb-534}
cancel each other. We shall use \ainf-functor identities
$b\phi=\phi b$, $b\psi=\psi b$ for terms \eqref{eq-bfpuwb-533} and
\eqref{eq-fpubwb-535}. Being a cocategory homomorphism, $\phi$
satisfies the identity
\[ \sum_{a+e=h} \phi_{a\alpha}\tens\phi_{e\gamma}
= \bigl[\Delta(\phi\tens\phi)\bigr]_{h;\alpha,\gamma}
= \phi_{h,\alpha+\gamma}\Delta_{\alpha+\gamma;\alpha,\gamma}
\]
for all non-negative integers $h$, where $\Delta$ is the cut
comultiplication. Similarly for $\psi$. Using this identity in
\eqref{eq-ffpuww1b1b-531} we get the equation to verify:
\begin{align*}
&- \sum_{x+q+z=k}^{v,w} (\phi_{xv}\tens r_q\tens\psi_{zw})
\sum_{\alpha+y+\beta=v+1+w}^{\alpha\le v,\beta\le w}
(1^{\tens\alpha}\tens b^\ca_y\tens1^{\tens\beta})b^\ca_{\alpha+1+\beta}
\\
&= \sum_{x+q+z=k}^{v,w,\alpha} (\phi_{xv}\tens r_q\tens\psi_{zw})
(b^\ca_{v\alpha}\tens1^{\tens1+w})b^\ca_{\alpha+1+w} \\
&+ \sum_{x+q+z=k}^{v,w,\beta} (\phi_{xv}\tens r_q\tens\psi_{zw})
(1^{\tens v+1}\tens b^\ca_{w\beta})b^\ca_{v+1+\beta}.
\end{align*}
It follows from the identity $b^2\pr_1=0:T^{v+1+w}s\ca\to s\ca$ valid
for arbitrary non-negative integers $v$, $w$, which we may rewrite like
this:
\begin{equation*}
\sum_{\alpha+y+\beta=v+1+w}^{\alpha\le v,\beta\le w}
(1^{\tens\alpha}\tens b^\ca_y\tens1^{\tens\beta})b^\ca_{\alpha+1+\beta}
+ \sum_\alpha (b^\ca_{v\alpha}\tens1^{\tens1+w})b^\ca_{\alpha+1+w}
+ \sum_\beta (1^{\tens v+1}\tens b^\ca_{w\beta})b^\ca_{v+1+\beta} = 0.
\end{equation*}
So the lemma is proved.
\end{proof}

The proposition follows.
\end{proof}

Let us consider now the question, when the discussed chain map is
null-homotopic.

\begin{corollary}\label{cor-chain-to-Ainf-FQ-A-null-homotopic}
Let $\phi,\psi:\cf\cq\to\ca$ be \ainf-functors. Let $P$ be a complex of
$\kk$\n-modules. Let $u:P\to sA_\infty(\cf\cq,\ca)(\phi,\psi)$ be a
chain map. The set (possibly empty) of homotopies
$h:P\to sA_\infty(\cf\cq,\ca)(\phi,\psi)$, $\deg h=-1$, such that
$u=dh+hB_1$ is in bijection with the set of data $(h',h_k)_{k>1}$,
consisting of
\begin{enumerate}
\item a homotopy $h':P\to sA_1(\cq,\ca)(\phi,\psi)$, $\deg h'=-1$, such
that $dh'+h'B_1=u'$, where $u'$ is given by \eqref{eq-u'-urestr2};

\item $\kk$\n-linear maps
\[ h_k: P \to \prod_{X,Y\in\Ob\cq}
\Com\bigl((s\cf\cq)^{\tens k}(X,Y),s\ca(X\phi,Y\psi)\bigr)
\]
of degree $-1$ for all $k>1$.
\end{enumerate}
The bijection maps $h$ to $(h',h_k)_{k>1}$, where $h_k=h\cdot\pr_k$ and
\[ h' = \bigl(P \rTTo^h sA_\infty(\cf\cq,\ca)(\phi,\psi) \rEpi^\restr
sA_1(\cf\cq,\ca)(\phi,\psi) \rEpi^\restr sA_1(\cq,\ca)(\phi,\psi)
\bigr).
\]
The inverse bijection can be recovered from the recurrent formula
\begin{multline*}
(-)^pb_k(ph_1) = pu_k -(pd)h_k - \sum_{a+q+c=k}^{\alpha,\beta}
(\phi_{a\alpha}\tens ph_q\tens\psi_{c\beta})b_{\alpha+1+\beta} \\
-(-)^p \sum_{a+q+c=k}^{a+c>0}
(1^{\tens a}\tens b_q\tens1^{\tens c})(ph_{a+1+c}):
(s\cf\cq)^{\tens k} \to s\ca,
\end{multline*}
where $k>1$, \(p\in P\), and \(\phi_{a\alpha}\), \(\psi_{c\beta}\) are
matrix elements of $\phi$, $\psi$.
\end{corollary}

\begin{proof}
We shall apply \propref{pro-chain-to-Ainf-FQ-A} to the complex
\(\Cone(\id:P\to P)\) instead of $P$. The graded $\kk$\n-module
\(\Cone(\id_P)=P\oplus P[1]\) is equipped with the differential
\((q,ps)d=(qd+p,-pds)\), \(p,q\in P\). The chain maps
\(\overline{u}:\Cone(\id_P)\to C\) to an arbitrary complex $C$ are in
bijection with pairs \((u:P\to C,h:P\to C)\), where \(u=dh+hd\) and
\(\deg h=-1\). The pair
\((u,h)=(\inj_1\overline{u},s\inj_2\overline{u})\) is assigned to
\(\overline{u}\), and the map \(\overline{u}:P\oplus P[1]\to C\),
\((q,ps)\mapsto qu+ph\) is assigned to a pair \((u,h)\). Indeed,
\(\overline{u}\) being chain map is equivalent to
\[ (q,ps)d\overline{u} = qdu +pu -pdh = qud +phd =(q,ps)\overline{u}d,
\]
that is, to conditions \(du=ud\), \(u=dh+hd\).

Thus, for a fixed chain map \(u:P\to C\) the set of homotopies
\(h:P\to C\), such that \(u=dh+hd\), is in bijection with the set of
chain maps \(\overline{u}:\Cone(\id_P)\to C\) such that
\(\inj_1\overline{u}=u:P\to C\). Applying this statement to
\(u:P\to C=sA_\infty(\cf\cq,\ca)(\phi,\psi)\) we find by
\propref{pro-chain-to-Ainf-FQ-A} that the set of homotopies
$h:P\to sA_\infty(\cf\cq,\ca)(\phi,\psi)$ such that $u=dh+hB_1$ is in
bijection with the set of data $(\overline{u}',\overline{u}_k)_{k>1}$,
such that
\begin{alignat*}2
\overline{u}' &: \Cone(\id_P) \to sA_1(\cq,\ca)(\phi,\psi) \qquad
\text{is a chain map},\quad & \inj_1\overline{u}' &= u', \\
\overline{u}_k &: \Cone(\id_P) \to \prod_{X,Y\in\Ob\cq}
\Com\bigl((s\cf\cq)^{\tens k}(X,Y),s\ca(X\phi,Y\psi)\bigr),
\quad \deg\overline{u}_k = 0,\quad & \inj_1\overline{u}_k &= u_k,
\end{alignat*}
therefore, in bijection with the set of data
$(h',h_k)_{k>1}=(s\inj_2\overline{u}',s\inj_2\overline{u}_k)_{k>1}$, as
stated in corollary.
\end{proof}

\subsection{Restriction as an $A_\infty$-functor.}
Let $\cq$ be a ($\fu$\n-small) differential graded $\kk$\n-quiver.
Denote by $\cf\cq$ the free \ainf-category generated by $\cq$. Let
$\ca$ be a ($\fu$\n-small) unital \ainf-category. There is the
restriction strict \ainf-functor
\[ \restr: A_{\infty}(\cf\cq,\ca) \to A_1(\cq,\ca), \qquad
(f:\cf\cq\to\ca) \mapsto (\overline{f}=(f_1\big|_\cq):\cq\to\ca).
\]
In fact, it is the composition of two strict \ainf-functors:
 $A_{\infty}(\cf\cq,\ca) \rTTo^{\restr_{\infty,1}}
 A_1(\cf\cq,\ca)\to A_1(\cq,\ca)$,
where the second comes from the full embedding
$\cq\hookrightarrow\cf\cq$. Its first component is
\begin{align}
\restr_1: sA_{\infty}(\cf\cq,\ca)(f,g) &\to
sA_1(\cq,\ca)(\overline{f},\overline{g}), \label{eq-restr-A(FQA)-A1(QA)} \\
r = (r_0,r_1,\dots,r_n,\dots) &\mapsto (r_0,r_1|_{\cq}) = \overline{r}.
\notag
\end{align}

\begin{theorem}\label{thm-restr-equivalence}
The \ainf-functor $\restr:A_\infty(\cf\cq,\ca)\to A_1(\cq,\ca)$ is an
equivalence.
\end{theorem}

\begin{proof}
Let us prove that restriction map \eqref{eq-restr-A(FQA)-A1(QA)}
is homotopy invertible. We construct a chain map going in the opposite
direction
\begin{equation*}
u: sA_1(\cq,\ca)(\overline{f},\overline{g}) \to
sA_\infty(\cf\cq,\ca)(f,g)
\end{equation*}
via \propref{pro-chain-to-Ainf-FQ-A} taking
$P=sA_1(\cq,\ca)(\overline{f},\overline{g})$. We choose
\begin{equation*}
u': sA_1(\cq,\ca)(\overline{f},\overline{g}) \to
sA_1(\cq,\ca)(\overline{f},\overline{g})
\end{equation*}
to be the identity map and $u_k=0$ for $k>1$. Therefore,
\[ u\cdot\restr_1 =u' = \id_{sA_1(\cq,\ca)(\overline{f},\overline{g})}.
\]
Denote
\[ v = \id_{sA_\infty(\cf\cq,\ca)(f,g)}
- \bigl[sA_\infty(\cf\cq,\ca)(f,g) \rTTo^{\restr_1}
sA_1(\cq,\ca)(\overline{f},\overline{g}) \rTTo^u
sA_\infty(\cf\cq,\ca)(f,g) \bigr].
\]
Let us prove that $v$ is null-homotopic via
\corref{cor-chain-to-Ainf-FQ-A-null-homotopic}, taking
$P=sA_\infty(\cf\cq,\ca)(f,g)$. A homotopy
$h:sA_\infty(\cf\cq,\ca)(f,g)\to sA_\infty(\cf\cq,\ca)(f,g)$,
$\deg h=-1$, such that $v=B_1h+hB_1$ is specified by
 $h'=0:sA_\infty(\cf\cq,\ca)(f,g)\to
 sA_1(\cq,\ca)(\overline{f},\overline{g})$
and $h_k=0$ for $k>1$. Indeed,
\[ v' = v\cdot\restr_1 = \restr_1 - \restr_1\cdot u\cdot\restr_1
= \restr_1 - \restr_1 = 0,
\]
so $v'=B_1h'+h'B_1$ and condition 1 of
\corref{cor-chain-to-Ainf-FQ-A-null-homotopic} is satisfied%
\footnote{By the way, the only non-vanishing component of $h$ is $h_1$.}.
Therefore, $u$ is homotopy inverse to $\restr_1$.

Let $\uni^\ca$ be a unit transformation of the unital \ainf-category
$\ca$. Then $A_1(\cq,\ca)$ is a unital \ainf-category with the unit
transformation $(1\tens\uni^\ca)M$ (cf.
\cite[Proposition~7.7]{Lyu-AinfCat}). The unit element for an object
$\phi\in\Ob A_1(\cq,\ca)$ is
$\sS{_\phi}\uni^{A_1(\cq,\ca)}_0:\kk\to sA_1(\cq,\ca)$,
$1\mapsto\phi\uni^\ca$. The \ainf-category $A_\infty(\cf\cq,\ca)$ is
also unital. To establish equivalence of these two \ainf-categories via
$\restr:A_\infty(\cf\cq,\ca)\to A_1(\cq,\ca)$ we verify the conditions
of Theorem~8.8 from \cite{Lyu-AinfCat}.

Consider the mapping $\Ob A_1(\cq,\ca)\to\Ob A_\infty(\cf\cq,\ca)$,
$\phi\mapsto\widehat{\phi}$, which extends a given chain map to a
strict \ainf-functor, constructed in \corref{cor-restr-strAinf}.
Clearly, $\overline{\widehat{\phi}}=\phi$. It remains to give two
mutually inverse cycles, which we choose as follows:
\begin{align*}
\sS{_\phi}r_0: \kk &\to sA_1(\cq,\ca)(\phi,\overline{\widehat{\phi}}),
\qquad 1 \mapsto \phi\uni^\ca, \\
\sS{_\phi}p_0: \kk &\to sA_1(\cq,\ca)(\overline{\widehat{\phi}},\phi),
\qquad 1 \mapsto \phi\uni^\ca.
\end{align*}
Clearly, $\sS{_\phi}r_0B_1=0$, $\sS{_\phi}p_0B_1=0$,
\begin{align*}
(\sS{_\phi}r_0\tens\sS{_\phi}p_0)B_2 - \sS{_\phi}\uni^{A_1(\cq,\ca)}_0:
1 &\mapsto (\phi\uni^\ca\tens\phi\uni^\ca)B_2 -\phi\uni^\ca \in\im B_1,
\\
(\sS{_\phi}p_0\tens\sS{_\phi}r_0)B_2 - \sS{_\phi}\uni^{A_1(\cq,\ca)}_0:
1 &\mapsto (\phi\uni^\ca\tens\phi\uni^\ca)B_2 -\phi\uni^\ca \in\im B_1.
\end{align*}
Therefore, all assumptions of Theorem~8.8 \cite{Lyu-AinfCat} are
satisfied. Thus, $\restr:A_\infty(\cf\cq,\ca)\to A_1(\cq,\ca)$ is an
\ainf-equivalence.
\end{proof}

\begin{corollary}
Every \ainf-functor $f:\cf\cq\to\ca$ is isomorphic to the strict
\ainf-functor $\wh{\overline{f}}:\cf\cq\to\ca$.
\end{corollary}

\begin{proof}
Note that $\overline{f}=\overline{\wh{\overline{f}}}$. The
$A_1$\n-transformation
 $\overline{f}\uni^\ca:\overline{f}\to\overline{\wh{\overline{f}}}
 :\cq\to\ca$
with the components $(\sS{_{Xf}}\uni^\ca_0,\overline{f}_1\uni^\ca_1)$
is natural. It is mapped by $u$ into a natural \ainf-transformation
$(\overline{f}\uni^\ca)u:f\to\wh{\overline{f}}:\cf\cq\to\ca$. Its zero
component $\sS{_{Xf}}\uni^\ca_0$ is invertible, therefore
$(\overline{f}\uni^\ca)u$ is invertible by
\cite[Proposition~7.15]{Lyu-AinfCat}.
\end{proof}

\section{\texorpdfstring{Representable 2-functors $A_\infty^u\to A_\infty^u$}
 {Representable 2-functors A-infinity-u to A-infinity-u}}
Recall that unital \ainf-categories, unital \ainf-functors and
equivalence classes of natural \ainf-transformations form a
2\n-category \cite{Lyu-AinfCat}. In order to distinguish between the
\ainf-category $A_\infty^u(\cc,\cd)$ and the ordinary category, whose
morphisms are equivalence classes of natural \ainf-transformations, we
denote the latter by
\[ \overline{A_\infty^u}(\cc,\cd) = H^0(A_\infty^u(\cc,\cd),m_1).
\]
The corresponding notation for the 2\n-category is
$\overline{A_\infty^u}$. We will see that arbitrary $A_N$\n-categories
can be viewed as 2\n-functors
$\overline{A_\infty^u}\to\overline{A_\infty^u}$. Moreover, they come
from certain generalizations called \ainfu-2-functors. There is a notion
of representability of such 2\n-functors, which explains some
constructions of \ainf-categories. For instance, a differential graded
$\kk$\n-quiver $\cq$ will be represented by the free \ainf-category
$\cf\cq$ generated by it.

\begin{definition}\label{def-weak-2-func}
A (strict) \emph{\ainfu-2-functor} $F:A_\infty^u\to A_\infty^u$ consists
of
\begin{enumerate}
\item a map $F:\Ob A_\infty^u\to\Ob A_\infty^u$;
\item a unital \ainf-functor
 $F=F_{\cc,\cd}:A_\infty^u(\cc,\cd)\to A_\infty^u(F\cc,F\cd)$
for each pair $\cc,\cd$ of unital \ainf-categories;

such that
\item $\id_{F\cc}=F(\id_\cc)$ for any unital \ainf-category $\cc$;
\item the equation
\begin{diagram}[LaTeXeqno]
TsA_\infty^u(\cc,\cd)\boxtimes TsA_\infty^u(\cd,\ce) & \rTTo^M &
TsA_\infty^u(\cc,\ce) \\
\dTTo<{F_{\cc,\cd}\boxtimes F_{\cd,\ce}} & = & \dTTo>{F_{\cc,\ce}} \\
TsA_\infty^u(F\cc,F\cd)\boxtimes TsA_\infty^u(F\cd,F\ce) & \rTTo^M &
TsA_\infty^u(F\cc,F\ce)
\label{dia-ainf-2-functor}
\end{diagram}
holds strictly for each triple $\cc,\cd,\ce$ of unital
\ainf-categories.
\end{enumerate}
\end{definition}

The \ainf-functor $F:A_\infty^u(\cc,\cd)\to A_\infty^u(F\cc,F\cd)$
consists of the mapping of objects
\[ \Ob F: \Ob A_\infty^u(\cc,\cd) \to \Ob A_\infty^u(F\cc,F\cd),
\qquad f \mapsto Ff,
\]
and the components $F_k$, $k\ge1$:
\begin{gather*}
F_1: sA_\infty^u(\cc,\cd)(f,g) \to sA_\infty^u(F\cc,F\cd)(Ff,Fg), \\
F_2: sA_\infty^u(\cc,\cd)(f,g)\tens sA_\infty^u(\cc,\cd)(g,h) \to
sA_\infty^u(F\cc,F\cd)(Ff,Fh),
\end{gather*}
and so on.

Weak versions of \ainfu-2-functors and 2-transformations between them
might be considered elsewhere.

\begin{definition}\label{def-weak-2-trans}
A (strict) \emph{\ainfu-2-transformation}
$\lambda:F\to G:A_\infty^u\to A_\infty^u$ of strict \ainfu-2-functors is
\begin{enumerate}
\item a family of unital \ainf-functors $\lambda_\cc:F\cc\to G\cc$,
$\cc\in\Ob A_\infty^u$;

such that
\item the diagram of \ainf-functors
\begin{diagram}[LaTeXeqno]
A_\infty^u(\cc,\cd) & \rTTo^F & A_\infty^u(F\cc,F\cd) \\
\dTTo<G & = & \dTTo>{(1\boxtimes\lambda_\cd)M} \\
A_\infty^u(G\cc,G\cd) & \rTTo^{(\lambda_\cc\boxtimes1)M} &
A_\infty^u(F\cc,G\cd)
\label{dia-ainf-2-transformation}
\end{diagram}
strictly commutes.
\end{enumerate}
An \ainfu-2-transformation $\lambda=(\lambda_\cc)$ for which
$\lambda_\cc$ are \ainf-equivalences is called a \emph{natural
\ainfu-2-equivalence}.
\end{definition}

Let us show now that the above notions induce ordinary strict
2-functors and strict 2-transformations in 0\n-th cohomology. Recall
that the strict 2\n-category $\overline{A_\infty^u}$ consists of
objects -- unital \ainf-categories, the category
$\overline{A_\infty^u}(\cc,\cd)$ for any pair of objects $\cc$, $\cd$,
the identity functor $\id_\cc$ for any unital \ainf-category $\cc$, and
the composition functor \cite{Lyu-AinfCat}
\begin{align*}
\overline{A_\infty^u}(\cc,\cd)(f,g)\times
\overline{A_\infty^u}(\cd,\ce)(h,k)
& \rTTo^{\bull^2} \overline{A_\infty^u}(\cc,\ce)(fh,gk), \\
(rs^{-1},ps^{-1}) &\rMapsTo (rhs^{-1}\tens gps^{-1})m_2.
\end{align*}
Given a strict \ainfu-2-functor $F$ as in \defref{def-weak-2-func} we
construct from it an ordinary strict 2\n-functor
 $\overline{F}=F:\Ob\overline{A_\infty^u}\to\Ob\overline{A_\infty^u}$,
 $\overline{F}=H^0(sF_1s^{-1}):\overline{A_\infty^u}(\cc,\cd)\to
 \overline{A_\infty^u}(F\cc,F\cd)$
as follows.

Denote
\begin{multline}
M_{10}\odot M_{01} = \bigl\{
sA_\infty^u(\cc,\cd)\boxtimes sA_\infty^u(\cd,\ce)
\rTTo^{\Delta_{10}\boxtimes\Delta_{01}}_\sim \\
[sA_\infty^u(\cc,\cd)\tens T^0sA_\infty^u(\cc,\cd)]\boxtimes
[T^0sA_\infty^u(\cd,\ce)\tens sA_\infty^u(\cd,\ce)] \\
\rTTo^\sim \relax
[sA_\infty^u(\cc,\cd)\boxtimes T^0sA_\infty^u(\cd,\ce)]\tens
[T^0sA_\infty^u(\cc,\cd)\boxtimes sA_\infty^u(\cd,\ce)] \\
\rTTo^{M_{10}\tens M_{01}}
sA_\infty^u(\cc,\ce)\tens sA_\infty^u(\cc,\ce) \bigr\},
\label{eq-M10-odot-M01}
\end{multline}
where the obvious isomorphisms $\Delta_{10}$ and $\Delta_{01}$ are
components of the comultiplication $\Delta$, the middle isomorphism is
that of distributivity law~\eqref{eq-Distributivity-law}, and the
components $M_{10}$ and $M_{01}$ of $M$ are the composition maps.

Property~\eqref{dia-ainf-2-functor} of $F$ implies that
\begin{equation}
(M_{10}\odot M_{01})(F_1\tens F_1)
= (F_1\boxtimes F_1)(M_{10}\odot M_{01}).
\label{eq-MMFF-FFMM}
\end{equation}
Indeed, the following diagram commutes
\begin{diagram}[nobalance]
sA_\infty^u(\cc,\cd)\boxtimes T^0sA_\infty^u(\cd,\ce)
& \rTTo^{M_{10}} & sA_\infty^u(\cc,\ce) \\
\dTTo<{F_1\boxtimes\Ob F} && \dTTo>{F_1} \\
sA_\infty^u(F\cc,F\cd)\boxtimes T^0sA_\infty^u(F\cd,F\ce)
& \rTTo^{M_{10}} & sA_\infty^u(F\cc,F\ce)
\end{diagram}
due to \eqref{dia-ainf-2-functor}. $\tens$\n-tensoring it with one
more similar diagram we get
\[ (M_{10}\tens M_{01})(F_1\tens F_1)
= [(F_1\boxtimes\Ob F)\tens(\Ob F\boxtimes F_1)](M_{10}\tens M_{01}).
\]
The isomorphisms in \eqref{eq-M10-odot-M01} commute with $F$ in
expected way, so \eqref{eq-MMFF-FFMM} follows.

We claim that the diagram
\begin{diagram}[LaTeXeqno]
sA_\infty^u(\cc,\cd)\boxtimes sA_\infty^u(\cd,\ce)
& \rTTo^{(M_{10}\odot M_{01})B_2} & sA_\infty^u(\cc,\ce) \\
\dTTo<{F_1\boxtimes F_1} && \dTTo>{F_1} \\
sA_\infty^u(F\cc,F\cd)\boxtimes sA_\infty^u(F\cd,F\ce)
& \rTTo^{(M_{10}\odot M_{01})B_2} & sA_\infty^u(F\cc,F\ce)
\label{dia-(MM)BF-FF(MM)B}
\end{diagram}
homotopically commutes. Indeed, since
\[ (1\tens B_1+B_1\tens1)F_2 + B_2F_1 = (F_1\tens F_1)B_2 + F_2B_1,
\]
we get
\begin{align*}
(&M_{10}\odot M_{01})B_2F_1 \\
&= (M_{10}\odot M_{01})(F_1\tens F_1)B_2
+ (M_{10}\odot M_{01})F_2B_1
- (M_{10}\odot M_{01})(1\tens B_1+B_1\tens1)F_2 \\
&= (F_1\boxtimes F_1)(M_{10}\odot M_{01})B_2
+ (M_{10}\odot M_{01})F_2B_1
- (1\boxtimes B_1+B_1\boxtimes1)(M_{10}\odot M_{01})F_2.
\end{align*}
We have used equations
\begin{align*}
(M_{10}\odot M_{01})(1\tens B_1)
&= (1\boxtimes B_1)(M_{10}\odot M_{01}), \\
(M_{10}\odot M_{01})(B_1\tens1) &= (B_1\boxtimes1)(M_{10}\odot M_{01}),
\end{align*}
which can be proved similarly to \eqref{eq-MMFF-FFMM} due to $M$ being
an \ainf-functor. Passing to cohomology we get from
\eqref{dia-(MM)BF-FF(MM)B} a strictly commutative diagram of functors
\begin{diagram}
H^0(A_\infty^u(\cc,\cd)\boxtimes A_\infty^u(\cd,\ce))
& \rTTo^{\bull^2} & H^0(A_\infty^u(\cc,\ce)) \\
\dTTo<{H^0(sF_1s^{-1}\boxtimes sF_1s^{-1})} &=& \dTTo>{H^0(sF_1s^{-1})}
\\
H^0(A_\infty^u(F\cc,F\cd)\boxtimes A_\infty^u(F\cd,F\ce))
& \rTTo^{\bull^2} & H^0(A_\infty^u(F\cc,F\ce))
\end{diagram}
since $\bull^2=H^0((s\boxtimes s)(M_{10}\odot M_{01})B_2s^{-1})$. Using
the K\"unneth map we come to strictly commutative diagram of functors
\begin{diagram}
\overline{A_\infty^u}(\cc,\cd)\times \overline{A_\infty^u}(\cd,\ce)
& \rTTo^{\bull^2} & \overline{A_\infty^u}(\cc,\ce) \\
\dTTo<{H^0(sF_1s^{-1})\times H^0(sF_1s^{-1})} & = &
\dTTo>{H^0(sF_1s^{-1})} \\
\overline{A_\infty^u}(F\cc,F\cd)\times \overline{A_\infty^u}(F\cd,F\ce)
& \rTTo^{\bull^2} & \overline{A_\infty^u}(F\cc,F\ce)
\end{diagram}
that is, to a usual strict 2\n-functor
$\overline{F}:\overline{A_\infty^u}\to \overline{A_\infty^u}$.

Let us show that an \ainfu-2-transformation
$\lambda:F\to G:A_\infty^u\to A_\infty^u$ as in
\defref{def-weak-2-trans} induces an ordinary strict 2\n-transformation
 $\overline{\lambda}:\overline{F}\to\overline{G}:
 \overline{A_\infty^u}\to\overline{A_\infty^u}$
in cohomology. Indeed, diagram~\eqref{dia-ainf-2-transformation}
implies commutativity of diagram
\begin{diagram}
sA_\infty^u(\cc,\cd) & \rTTo^{F_1} & sA_\infty^u(F\cc,F\cd) \\
\dTTo<{G_1} & = & \dTTo>{(1\boxtimes\lambda_\cd)M_{10}} \\
sA_\infty^u(G\cc,G\cd) & \rTTo^{(\lambda_\cc\boxtimes1)M} &
sA_\infty^u(F\cc,G\cd) \\
r:f\to g:G\cc\to G\cd & \rMapsTo &
\lambda_\cc r:\lambda_\cc f\to\lambda_\cc g:F\cc\to G\cd.
\end{diagram}
Passing to cohomology we get
\begin{diagram}
\overline{A_\infty^u}(\cc,\cd) & \rTTo^{H^0(sF_1s^{-1})}
& \overline{A_\infty^u}(F\cc,F\cd) \\
\dTTo<{H^0(sG_1s^{-1})} & = &
\dTTo<{\_\cdot\lambda_\cd}>{=\overline{A_\infty^u}(F\cc,\lambda_\cd)}\\
\overline{A_\infty^u}(G\cc,G\cd) &
\rTTo^{\lambda_\cc\cdot\_}_{\overline{A_\infty^u}(\lambda_\cc,G\cd)} &
\overline{A_\infty^u}(F\cc,G\cd)
\end{diagram}
Therefore,
$\overline{\lambda_\cc}\in\Ob\overline{A_\infty^u}(F\cc,G\cc)$ form a
strict 2\n-transformation
 $\overline{\lambda}:\overline{F}\to\overline{G}:
 \overline{A_\infty^u}\to\overline{A_\infty^u}$.

\subsection{Examples of $A_\infty^u$-2-functors.}
Let $\ca$ be an $A_N$\n-category, $1\le N\le\infty$. It determines an
\ainfu-2-functor $F=A_N(\ca,\_):A_{\infty}^u\to A_{\infty}^u$, given by
the following data:
\begin{enumerate}
\item the map $F:\Ob A_\infty^u\to\Ob A_\infty^u$,
$\cc\mapsto A_N(\ca,\cc)$ (the category $A_N(\ca,\cc)$ is unital by
\cite[Proposition~7.7]{Lyu-AinfCat});

\item the unital strict \ainf-functor
 $F=A_N(\ca,\_):A_\infty^u(\cc,\cd)\to
 A_\infty^u(A_N(\ca,\cc),A_N(\ca,\cd))$
for each pair $\cc,\cd$ of unital \ainf-categories (cf.
\cite[Propositions 6.2, 8.4]{Lyu-AinfCat}).
\end{enumerate}

Clearly, $\id_{A_N(\ca,\cc)}=(1\boxtimes\id_{\cc})M=A_N(\ca,\id_\cc)$. We
want to prove now that the equation
\begin{multline}
\bigl[ TsA_\infty^u(\cc,\cd)\boxtimes TsA_\infty^u(\cd,\ce) \rTTo^M
TsA_\infty^u(\cc,\ce) \rTTo^{A_N(\ca,\_)}
TsA_\infty^u(A_N(\ca,\cc),A_N(\ca,\ce)) \bigr] \\
= \bigl[ TsA_\infty^u(\cc,\cd)\boxtimes TsA_\infty^u(\cd,\ce)
\rTTo^{A_N(\ca,\_)\boxtimes A_N(\ca,\_)} \\
TsA_\infty^u(A_N(\ca,\cc),A_N(\ca,\cd))\boxtimes
TsA_\infty^u(A_N(\ca,\cd),A_N(\ca,\ce)) \\
\rTTo^M TsA_\infty^u(A_N(\ca,\cc),A_N(\ca,\ce)) \bigr]
\label{eq-7TsA}
\end{multline}
holds strictly for each triple $\cc,\cd,\ce$ of unital
\ainf-categories. In fact, this $F$ is a restriction of an
\ainf-2-functor $F:\Ob A_\infty\to\Ob A_\infty$,
$\cc\mapsto A_N(\ca,\cc)$, which is defined just as in
\defref{def-weak-2-func} without mentioning the unitality.
Equation~\eqref{eq-7TsA} follows from a similar equation without the
unitality index $u$. To prove it we consider the compositions
\begin{multline*}
\bigl[ TsA_N(\ca,\cc)\boxtimes TsA_\infty(\cc,\cd)\boxtimes TsA_\infty(\cd,\ce)
\rTTo^{1\boxtimes M} TsA_N(\ca,\cc)\boxtimes TsA_\infty(\cc,\ce) \\
\rTTo^{1\boxtimes A_N(\ca,\_)}
TsA_N(\ca,\cc)\boxtimes TsA_\infty(A_N(\ca,\cc),A_N(\ca,\ce)) \rTTo^\alpha
TsA_N(\ca,\ce) \bigr] \\
= \bigl[
TsA_N(\ca,\cc)\boxtimes TsA_\infty(\cc,\cd)\boxtimes TsA_\infty(\cd,\ce) \\
\rTTo^{1\boxtimes M} TsA_N(\ca,\cc)\boxtimes TsA_\infty(\cc,\ce) \rTTo^M
TsA_N(\ca,\ce) \bigr] \\
= \bigl[
TsA_N(\ca,\cc)\boxtimes TsA_\infty(\cc,\cd)\boxtimes TsA_\infty(\cd,\ce) \\
\rTTo^{M\boxtimes1} TsA_N(\ca,\cc)\boxtimes TsA_\infty(\cc,\ce) \rTTo^M
TsA_N(\ca,\ce) \bigr] \\
= \bigl[
TsA_N(\ca,\cc)\boxtimes TsA_\infty(\cc,\cd)\boxtimes TsA_\infty(\cd,\ce)
\rTTo^{1\boxtimes A_N(\ca,\_)\boxtimes1} \\
TsA_N(\ca,\cc)\boxtimes TsA_\infty(A_N(\ca,\cc),A_N(\ca,\cd))
\boxtimes TsA_\infty(\cd,\ce) \\
\rTTo^{\alpha\boxtimes1} TsA_N(\ca,\cd)\boxtimes TsA_\infty(\cd,\ce) \\
\rTTo^{1\boxtimes A_N(\ca,\_)}
TsA_N(\ca,\cd)\boxtimes TsA_\infty(A_N(\ca,\cd),A_N(\ca,\ce))
\rTTo^\alpha TsA_N(\ca,\ce) \bigr] \\
= \bigl[
TsA_N(\ca,\cc)\boxtimes TsA_\infty(\cc,\cd)\boxtimes TsA_\infty(\cd,\ce)
\rTTo^{1\boxtimes A_N(\ca,\_)\boxtimes A_N(\ca,\_)} \\
TsA_N(\ca,\cc)\boxtimes TsA_\infty(A_N(\ca,\cc),A_N(\ca,\cd))
\boxtimes TsA_\infty(A_N(\ca,\cd),A_N(\ca,\ce)) \\
\rTTo^{1\boxtimes M}
TsA_N(\ca,\cc)\boxtimes TsA_\infty(A_N(\ca,\cc),A_N(\ca,\ce))
\rTTo^\alpha TsA_N(\ca,\ce) \bigr].
\end{multline*}
By \propref{prop-phi-C1-Cq-psi-ainf} we deduce equation~\eqref{eq-7TsA}
(see also \cite[Proposition~5.5]{Lyu-AinfCat}).

Let now $\ca$ be a unital \ainf-category. It determines an
\ainfu-2-functor $G=A_\infty^u(\ca,\_):A_\infty^u\to A_\infty^u$, given
by the following data:
\begin{enumerate}
\item the map $G:\Ob A_\infty^u\to\Ob A_\infty^u$,
$\cc\mapsto A_\infty^u(\ca,\cc)$ (the category $A_\infty^u(\ca,\cc)$ is
unital by \cite[Proposition~7.7]{Lyu-AinfCat});

\item the unital strict \ainf-functor
 $G=A_\infty^u(\ca,\_):A_\infty^u(\cc,\cd)\to
 A_\infty^u(A_\infty^u(\ca,\cc),A_\infty^u(\ca,\cd))$
for each pair $\cc,\cd$ of unital \ainf-categories, determined from
\begin{multline*}
M = \bigl[ TsA_\infty^u(\ca,\cb)\boxtimes TsA_\infty^u(\cb,\cc)
\rTTo^{1\boxtimes A_\infty^u(\ca,\_)} \\
TsA_\infty^u(\ca,\cb)\boxtimes
TsA_\infty^u(A_\infty^u(\ca,\cb),A_\infty^u(\ca,\cc))
\rTTo^\alpha TsA_\infty^u(\ca,\cc)\bigr].
\end{multline*}
(cf. \cite[Propositions 6.2, 8.4]{Lyu-AinfCat}).
\end{enumerate}
Clearly, $G\cc$ are full \ainf-subcategories of $F\cc$ for the
\ainfu-2-functor $F=A_\infty(\ca,\_)$. Furthermore, \ainf-functors
$G_{\cc,\cd}(f)$ are restrictions of \ainf-functors $F_{\cc,\cd}(f)$,
so $G$ is a full \ainfu-2-subfunctor of $F$. In particular, $G$
satisfies equation~\eqref{dia-ainf-2-functor}. Another way to prove
that $G$ is an \ainfu-2-functor is to repeat the reasoning concerning
$F$.

\subsection{Example of an $A_\infty^u$-2-equivalence.}
Assume that $\cq$ is a differential graded $\kk$\n-quiver. As usual,
$\cf\cq$ denotes the free \ainf-category generated by it. We claim that
 $\restr:A_\infty(\cf\cq,\_)\to A_1(\cq,\_):A_\infty^u\to A_\infty^u$
is a strict 2\n-natural \ainf-equivalence. Indeed, it is given by the
family of unital \ainf-functors
 $\restr_\cc:A_\infty(\cf\cq,\cc)\to A_1(\cq,\cc)$,
$\cc\in\Ob A_\infty^u$, which are equivalences by
\thmref{thm-restr-equivalence}. We have to prove that the diagram of
\ainf-functors
\begin{diagram}[LaTeXeqno]
A_\infty^u(\cc,\cd) & \rTTo^{A_\infty(\cf\cq,\_)} &
A_\infty^u(A_\infty(\cf\cq,\cc),A_\infty(\cf\cq,\cd)) \\
\dTTo<{A_1(\cq,\_)} & = & \dTTo>{(1\boxtimes\restr_\cd)M} \\
A_\infty^u(A_1(\cq,\cc),A_1(\cq,\cd)) & \rTTo^{(\restr_\cc\boxtimes1)M} &
A_\infty^u(A_\infty(\cf\cq,\cc),A_1(\cq,\cd))
\label{dia-Ainf-A1-rest-rest}
\end{diagram}
commutes. Notice that all arrows in this diagram are strict
\ainf-functors. Indeed, $A_\infty(\cf\cq,\_)$ and $A_1(\cq,\_)$ are
strict by \cite[Proposition~6.2]{Lyu-AinfCat}. For an arbitrary
\ainf-functor $f$ the components $[(f\boxtimes1)M]_n=(f\boxtimes1)M_{0n}$
vanish for all $n$ except for $n=1$, thus, $(f\boxtimes1)M$ is strict. The
\ainf-functor $g=\restr_\cd$ is strict, hence, the $n$\n-th component
\[ [(1\boxtimes g)M]_n: r^1\tdt r^n \mapsto (r^1\tdt r^n \mid g)M_{n0}
\]
of the \ainf-functor $(1\boxtimes g)M$ satisfies the equation
\[ [(r^1\tdt r^n \mid g)M_{n0}]_k = (r^1\tdt r^n)\theta_{k1}g_1.
\]
If the right hand side does not vanish, then $n\le1\le k+n$, so $n=1$
and $(1\boxtimes g)M$ is strict.

Given an \ainf-transformation $t:g\to h:\cc\to\cd$ between unital
\ainf-functors we find that
\begin{align*}
A_\infty(\cf\cq,\_)(t) = [(1\boxtimes t)M &:(1\boxtimes g)M\to(1\boxtimes h)M:
A_\infty(\cf\cq,\cc)\to A_\infty(\cf\cq,\cd)], \\
A_1(\cq,\_)(t) = [(1\boxtimes t)M &:(1\boxtimes g)M\to(1\boxtimes h)M:
A_1(\cq,\cc)\to A_1(\cq,\cd)], \\
[(1\boxtimes\restr_\cd)M]A_\infty(\cf\cq,\_)(t) &=
[((1\boxtimes t)M)\cdot\restr_\cd: ((1\boxtimes g)M)\cdot\restr_\cd \\
&\hspace*{3em} \to((1\boxtimes h)M)\cdot\restr_\cd:
A_\infty(\cf\cq,\cc)\to A_1(\cq,\cd)], \\
[(\restr_\cc\boxtimes1)M]A_1(\cq,\_)(t) &=
[\restr_\cc\cdot((1\boxtimes t)M):\restr_\cc\cdot((1\boxtimes g)M) \\
&\hspace*{3em} \to\restr_\cc\cdot((1\boxtimes h)M):
A_\infty(\cf\cq,\cc)\to A_1(\cq,\cd)].
\end{align*}
We have to verify that the last two \ainf-transformations are equal.
First of all, let us show that mappings of objects in
\eqref{dia-Ainf-A1-rest-rest} commute. Given a unital \ainf-functor
$g:\cc\to\cd$, we are going to check that
\begin{equation}
[(1\boxtimes g)M]_n\cdot\restr_1=\restr^{\tens n}_1\cdot[(1\boxtimes g)M]_n
\label{eq-1gMn-rest-rest-1gMn}
\end{equation}
for any $n\ge1$. Indeed, for any $n$\n-tuple of composable
\ainf-transformations
\[ f^0 \rTTo^{r^1} f^1 \rTTo \dots \rTTo^{r^n} f^n:\cf\cq\to\cc,
\]
we have in both cases
\begin{align*}
\{(r^1\tens\dots\tens r^n)[(1\boxtimes g)M]_n\}_0 &=
[(r^1\tens\dots\tens r^n|g)M_{n0}]_0=(r^1_0\tens\dots\tens r^n_0)g_n,\\
\{(r^1\tens\dots\tens r^n)[(1\boxtimes g)M]_n\}_1 &=
[(r^1\tens\dots\tens r^n|g)M_{n0}]_1 \\
&= \sum_{i=1}^n(r^1_0\tens\dots\tens r^{i-1}_0\tens
r^i_1\tens r^{i+1}_0\tens\dots\tens r^n_0)g_n \\
&+ \sum_{i=1}^n(r^1_0\tens\dots\tens r^{i-1}_0\tens f^i_1\tens
r^{i+1}_0\tens\dots\tens r^n_0)g_{n+1}.
\end{align*}
Note that the right hand sides depend only on 0\n-th and 1\n-st
components of $r^i$, $f^i$. This is precisely what is claimed by
equation~\eqref{eq-1gMn-rest-rest-1gMn}.

The coincidence of \ainf-transformations
$((1\boxtimes t)M)\cdot\restr_\cd=\restr_\cc\cdot((1\boxtimes t)M)$ follows
similarly from the computation:
\begin{align*}
\{(r^1\tens\dots\tens r^n)[(1\boxtimes t)M]_n\}_0 &=
[(r^1\tens\dots\tens r^n\boxtimes t)M_{n1}]_0
= (r^1_0\tens\dots\tens r^n_0)t_n, \\
\{(r^1\tens\dots\tens r^n)[(1\boxtimes t)M]_n\}_1 &=
[(r^1\tens\dots\tens r^n\boxtimes t)M_{n1}]_1 \\
&= \sum_{i=1}^n(r^1_0\tens\dots\tens r^{i-1}_0\tens
r^i_1\tens r^{i+1}_0\tens\dots\tens r^n_0)t_n \\
&+ \sum_{i=1}^n(r^1_0\tens\dots\tens r^{i-1}_0\tens f^i_1\tens
r^{i+1}_0\tens\dots\tens r^n_0)t_{n+1}.
\end{align*}

\subsection{Representability.}
An \ainfu-2-functor $F:A_\infty^u\to A_\infty^u$ is called
\emph{representable}, if it is naturally \ainfu-2-equivalent to the
\ainfu-2-functor $A_\infty(\ca,\_):A_\infty^u\to A_\infty^u$ for some
\ainf-category $\ca$. The above results imply that the \ainfu-2-functor
$A_1(\cq,\_)$ corresponding to a differential graded $\kk$\n-quiver
$\cq$ is represented by the free \ainf-category $\cf\cq$ generated by
$\cq$.

This definition of representability has a disadvantage: many different
\ainf-categories can represent the same \ainfu-2-functor. More
attractive notion is the following. An \ainfu-2-functor
$F:A_\infty^u\to A_\infty^u$ is called \emph{unitally representable},
if it is naturally \ainfu-2-equivalent to the \ainfu-2-functor
$A_\infty^u(\ca,\_):A_\infty^u\to A_\infty^u$ for some unital
\ainf-category $\ca$. Such $\ca$ is unique up to an \ainf-equivalence.
Indeed, composing a natural 2\n-equivalence
 $\overline{\lambda}:
 \overline{A_\infty^u(\ca,\_)}\to\overline{A_\infty^u(\cb,\_)}:
 \overline{A_\infty^u}\to\overline{A_\infty^u}$
with the 0\n-th cohomology 2\n-functor
$H^0:\overline{A_\infty^u}\to\Cat$, we get a natural 2\n-equivalence
 $H^0\overline{\lambda}:H^0\overline{A_\infty^u(\ca,\_)}\to
 H^0\overline{A_\infty^u(\cb,\_)}:\overline{A_\infty^u}\to\Cat$.
However,
$H^0\overline{A_\infty^u(\ca,\_)}=\overline{A_\infty^u}(\ca,\_)$, so
using a 2\n-category version of Yoneda lemma one can deduce that $\ca$
and $\cb$ are equivalent in $\overline{A_\infty^u}$. We shall present
an example of unital representability in subsequent
publication~\cite{math.CT/0306018}.

\subsection{Acknowledgements.}
We are grateful to all the participants of the \ainf-category seminar
at the Institute of Mathematics, Kyiv, for attention and fruitful
discussions, especially to Yu.~Bespalov and S.~Ovsienko. One of us
(V.L.) is grateful to Max-Planck-Institut f\"ur Mathematik for warm
hospitality and support at the final stage of this research.

\bibliographystyle{amsalpha}

\begin{thebibliography}{Kon95}

\bibitem[Fuk93]{Fukaya:A-infty}
K.~Fukaya, \emph{Morse homotopy, ${A}_\infty$-category, and {F}loer
  homologies}, {Proc. of GARC Workshop on Geometry and Topology '93}
  (H.~J. Kim, ed.), Lecture Notes, no.~18, Seoul Nat. Univ., Seoul,
  1993, pp.~1--102.

\bibitem[Fuk02]{Fukaya:FloerMirror-II}
K.~Fukaya, \emph{Floer homology and mirror symmetry. {II}}, {Minimal
  surfaces, geometric analysis and symplectic geometry (Baltimore, MD,
  1999)}, Adv. Stud. Pure Math., vol.~34, Math. Soc. Japan, Tokyo,
  2002, pp.~31--127.

\bibitem[Kad82]{Kadeishvili82}
T.~V. Kadeishvili, \emph{The algebraic structure in the homology of an
  ${A}(\infty)$-algebra}, Soobshch. Akad. Nauk Gruzin. SSR
  \textbf{108} (1982), no.~2, 249--252, in Russian.

\bibitem[Kel01]{math.RA/9910179}
B.~Keller, \emph{Introduction to {A}-infinity algebras and modules},
  Homology, Homotopy and Applications \textbf{3} (2001), no.~1, 1--35,
  \href{http://arXiv.org/abs/math.RA/9910179}{{\tt
  arXiv:\linebreak[1]math.RA/\linebreak[1]9910179}} ,
  \url{http://intlpress.com/HHA/v3/n1/a1/}.

\bibitem[Kon95]{Kontsevich:alg-geom/9411018}
M.~Kontsevich, \emph{Homological algebra of mirror symmetry}, Proc.
  Internat. Cong. Math., Z\"urich, Switzerland 1994 (Basel), vol.~1,
  Birkh\"auser Verlag, 1995, 120--139.

\bibitem[KS02]{KonSoi-AinfCat-NCgeom}
M.~Kontsevich and Y.~S. Soibelman, \emph{${A}_\infty$-categories and
  non-commutative geometry}, in preparation, 2002.

\bibitem[KS]{KonSoi-book}
M.~Kontsevich and Y.~S. Soibelman, \emph{Deformation theory}, book in
  preparation.

\bibitem[LH03]{Lefevre-Ainfty-these}
K.~Lef\`evre-Hasegawa, \emph{Sur les ${A}_\infty$-cat\'egories}, Ph.D.
  thesis, Universit\'e Paris 7, U.F.R. de Math\'ematiques, 2003,
  \href{http://arXiv.org/abs/math.CT/0310337}{{\tt
  arXiv:\linebreak[1]math.CT/\linebreak[1]0310337}}.

\bibitem[LM04]{math.CT/0306018}
V.~V. Lyubashenko and O.~Manzyuk, \emph{Quotients of unital
  ${A}_\infty$-categories}, 2004,
  \href{http://arXiv.org/abs/math.CT/0306018}{{\tt
  arXiv:\linebreak[1]math.CT/\linebreak[1]0306018}}.

\bibitem[LO02]{LyuOvs-iResAiFn}
V.~V. Lyubashenko and S.~A. Ovsienko, \emph{A construction of
  quotient $A_\infty$-categories},
  Homology, Homotopy Appl. \textbf{8} (2006), no.~2, 157--203,
  \href{http://arXiv.org/abs/math.CT/0211037}{{\tt
  arXiv:\linebreak[1]math.CT/\linebreak[1]0211037}}
  \url{http://intlpress.com/HHA/v8/n2/a9/}.

\bibitem[Lyu03]{Lyu-AinfCat}
V.~V. Lyubashenko, \emph{Category of ${A}_\infty$-categories},
  Homology, Homotopy and Applications \textbf{5} (2003), no.~1, 1--48,
  \href{http://arXiv.org/abs/math.CT/0210047}{{\tt
  arXiv:\linebreak[1]math.CT/\linebreak[1]0210047}} ,
  \url{http://intlpress.com/HHA/v5/n1/a1/}.

\bibitem[Sta63]{Stasheff:HomAssoc}
J.~D. Stasheff, \emph{Homotopy associativity of {H}-spaces, I $\&$ II},
  Trans. Amer. Math. Soc. \textbf{108} (1963), 275--292, 293--312.

\end{thebibliography}
\ifx\chooseClass2
\else
\tableofcontents
\fi
\end{document}